\documentclass[11pt]{amsart}
\usepackage{amsmath,amsfonts,amsthm,amsopn,amssymb,mathtools,stmaryrd}
\usepackage{cite,marginnote}
\usepackage{pdfsync}

\usepackage{color,enumitem,graphicx}
\usepackage[colorlinks=true,urlcolor=blue,
citecolor=red,linkcolor=blue,linktocpage,pdfpagelabels,
bookmarksnumbered,bookmarksopen]{hyperref}
\usepackage[english]{babel}

\usepackage[left=2.75cm,right=2.75cm,top=3cm,bottom=3cm]{geometry}
\usepackage[]{xcolor}

\numberwithin{equation}{section}

\makeindex
\newtheorem{theorem}{Theorem}[section]
\newtheorem{proposition}[theorem]{Proposition}
\newtheorem{lemma}[theorem]{Lemma}
\newtheorem{remark}[theorem]{Remark}
\newtheorem{example}[theorem]{Example}
\newtheorem{corollary}[theorem]{Corollary}
\newtheorem{definition}[theorem]{Definition}

\newcommand{\bt}{\begin{theorem}}
\newcommand{\et}{\end{theorem}}
\newcommand{\bl}{\begin{lemma}}
\newcommand{\el}{\end{lemma}}
\newcommand{\bd}{\begin{definition}}
\newcommand{\ed}{\end{definition}}
\newcommand{\bc}{\begin{corollary}}
\newcommand{\ec}{\end{corollary}}
\newcommand{\bp}{\begin{proof}}
\newcommand{\ep}{\end{proof}}
\newcommand{\bx}{\begin{example}}
\newcommand{\ex}{\end{example}}
\newcommand{\bi}{\begin{exercise}}
\newcommand{\ei}{\end{exercise}}
\newcommand{\bo}{\begin{proposition}}
\newcommand{\eo}{\end{proposition}}
\newcommand{\br}{\begin{remark}}
\newcommand{\er}{\end{remark}}
\newcommand{\be}{\begin{equation}}
\newcommand{\ee}{\end{equation}}
\newcommand{\ba}{\begin{align}}
\newcommand{\ea}{\end{align}}
\newcommand{\bn}{\begin{enumerate}}
\newcommand{\en}{\end{enumerate}}
\newcommand{\bg}{\begin{align*}}
\newcommand{\eg}{\end{align*}}
\newcommand{\bcs}{\begin{cases}}
\newcommand{\ecs}{\end{cases}}
\newcommand{\bean}{\begin{eqnarray*}}
\newcommand{\eean}{\end{eqnarray*}}

\newcommand{\s}{\section}

\newcommand{\na}{\nabla}

\newcommand{\R}{\mathbb R}

\newcommand{\e}{\varepsilon}

\title[Semiclassical states for coupled nonlinear Schr\"{o}dinger equations]{Semiclassical states for coupled nonlinear Schr\"{o}dinger equations with critical frequency}

\author[T. Y.\ Chen]{Taiyong Chen}
\author[Y. H.\ Jiang]{Yahui Jiang}
\author[M. Squassina]{Marco Squassina}
\author[J. J.\ Zhang]{Jianjun Zhang}

\address[T. Y.\ Chen]{\newline\indent School of Mathematics
\newline\indent
China University of Mining and Technology
\newline\indent
Xuzhou, 221116, China}
\email{\href{mailto:taiyongchencumt@163.com}{taiyongchencumt@163.com}}

\address[Y. H.\ Jiang]{\newline\indent School of Mathematics
\newline\indent
China University of Mining and Technology
\newline\indent
Xuzhou, 221116, China}
\email{\href{mailto:18843111149@163.com}{18843111149@163.com}}

\address[M.\ Squassina]{\newline\indent College of Science
\newline\indent
Princess Nourah Bint Abdul Rahman University
\newline\indent
Saudi Arabia, Riyadh, PO Box 84428}
\email{\href{mailto:marsquassina@pnu.edu.sa}{marsquassina@pnu.edu.sa}}

\address{Dipartimento di Matematica e Fisica
	\newline\indent
	Universit\`a Cattolica del Sacro Cuore
	\newline\indent
	Via dei Musei 41, Brescia, Italy}
\email{\href{mailto:marco.squassina@unicatt.it}{marco.squassina@unicatt.it}}

\address[J. J.\ Zhang]{\newline\indent College of Mathematica and Statistics
\newline\indent
Chongqing Jiaotong University
\newline\indent
Chongqing 400074, China}
\email{\href{mailto:zhangjianjun09@tsinghua.org.cn}{zhangjianjun09@tsinghua.org.cn}}

\thanks{Corresponding author: \texttt{marsquassina@pnu.edu.sa}}
\subjclass[2000]{35B05 35J50}
\date{}
\keywords{Nonlinear Schr\"{o}dinger systems, semiclassical limit, critical frequency, variational methods.}

\begin{document}

\begin{abstract}
In this paper, we are concerned with the coupled nonlinear Schr\"{o}dinger system
\begin{align*}
\begin{cases}
-\varepsilon^{2}\Delta u+a(x)u=\mu_{1}u^{3}+\beta v^{2}u \ \ \ \ \mbox{in}\ \mathbb{R}^{N},\\
-\varepsilon^{2}\Delta v+b(x)v=\mu_{2}v^{3}+\beta u^{2}v \ \ \ \ \ \mbox{in}\ \mathbb{R}^{N},
\end{cases}
\end{align*}
where $1\leq N\leq3$, $\mu_{1},\mu_{2},\beta>0$, $a(x)$ and $b(x)$ are nonnegative continuous potentials, and $\varepsilon>0$ is a small parameter. We show the existence of positive ground state solutions for the system above and also establish the concentration behaviour as $\varepsilon\rightarrow0$, when $a(x)$ and $b(x)$ achieve 0 with a homogeneous behaviour or vanish in some nonempty open set with smooth boundary.
\end{abstract}

\maketitle

\s{Introduction}
\renewcommand{\theequation}{1.\arabic{equation}}
\noindent {\bf 1.1 Background.}
Consider the following two-component coupled nonlinear Schr\"{o}dinger equations (also known as Gross-Pitaevskii equations):
\begin{align}
\label{question1}
\begin{cases}
-i\hbar\frac{\partial}{\partial_{t}}\Phi_{1}+V_{1}(x)\Phi_{1}=\frac{\hbar^{2}}{2m}\Delta\Phi_{1}+\mu_{1}|\Phi_{1}|^{2}\Phi_{1}+\beta|\Phi_{2}|^{2}\Phi_{1},\ \ \ \mbox{in}\ \mathbb{R}^{N},\ t>0,\\
-i\hbar\frac{\partial}{\partial_{t}}\Phi_{2}+V_{2}(x)\Phi_{2}=\frac{\hbar^{2}}{2m}\Delta\Phi_{2}+\mu_{2}|\Phi_{2}|^{2}\Phi_{2}+\beta|\Phi_{1}|^{2}\Phi_{2},\ \ \ \mbox{in}\ \mathbb{R}^{N},\ t>0,
\end{cases}
\end{align}
where $N\leq3$, $i$ is the imaginary unit, $\hbar$ is the Plank constant, $\mu_{1},\mu_{2}>0$ and $\beta\neq0$ is a coupling constant. System (\ref{question1}) has applications in many physical problems, such as nonlinear optics, continuum mechanics, plasma physics and Bose-Einstein condensates theory for multispecies Bose-Einstein condensates. In the context of Bose-Einstein condensates, $\Phi_{j}(x,t)(j=1,2)$ are the corresponding condensate amplitudes; $V_{j}(x)$ are the traping potentials; $\mu_{j}$ and $\beta$ are the intraspecies and interspecies scattering lengths. The physically realistic spatial dimension are $1\leq N\leq3$. When $N=2$, problem (\ref{question1}) arises in the Hartree-Fock theory for a double condenstate, $i.e.$, a binary mixture of Bose-Einstein condensates in two different hyperfine states $|1\rangle$ and $|2\rangle$ (see\cite{BCJJ,E}). The sign of $\beta$ determines whether the interactions of states are repulsive or attractive. The interaction is attractive if $\beta>0$ and repulsive if $\beta<0$, where the two states are in strong competition.\

A standing wave solution of (\ref{question1}) is of the following form
\begin{align}
\label{standing}
\Phi_{1}(x,t)=e^{-iE_{t}/\hbar}u(x)\ \ \mbox{and}\ \ \Phi_{2}(x,t)=e^{-iE_{t}/\hbar}v(x).
\end{align}
Substituting (\ref{standing}) into (\ref{question1}), then system \eqref{question1} is reduced to the following elliptic system
\begin{align}
\begin{cases}
-\frac{\hbar^{2}}{2m}\Delta u+(V_{1}(x)-E)u=\mu_{1}u^{3}+\beta v^{2}u,\ \ \ \mbox{in}\ \mathbb{R}^{N},\\
-\frac{\hbar^{2}}{2m}\Delta v+(V_{2}(x)-E)v=\mu_{2}v^{3}+\beta u^{2}v,\ \ \ \ \mbox{in}\ \mathbb{R}^{N}.\\
\end{cases}
\end{align}
Renaming the parameters by
\begin{align*}
\varepsilon=\sqrt{\frac{\hbar^{2}}{2m}},\ a(x)=V_{1}(x)-E,\ b(x)=V_{2}(x)-E,
\end{align*}
we obtain the folowing elliptic system
\begin{align*}
\label{question2}
\begin{cases}
-\varepsilon^{2}\Delta u+a(x)u=\mu_{1}u^{3}+\beta v^{2}u \ \ \ \ \mbox{in}\ \mathbb{R}^{N},\\
-\varepsilon^{2}\Delta v+b(x)v=\mu_{2}v^{3}+\beta u^{2}v \ \ \ \ \ \mbox{in}\ \mathbb{R}^{N}.\tag{$\mathcal{C}_{\varepsilon}$}
\end{cases}
\end{align*}
Before going further, we point out that the system (\ref{question2}) also possesses a trivial solution $(0,0)$ and semi-trivial solutions of type $(u,0)$ or $(0,v)$. A solution $(u,v)$ of (\ref{question2}) is nontrivial if $u\not\equiv0$ and $v\not\equiv0$. A solution $(u,v)$ with $u>0$ and $v>0$ is called a positive solution.\

When $\varepsilon=1$, $a(x)$ and $b(x)$ are positive constants, the existence, multiplicity, bifurcation and concentration behavior of positive solutions of (\ref{question2}) have been entensively studied, see Refs.\cite{AE,AM,TNZ,DO,N,TZ,Lin,LEB,Liu,Sirakov,WT}. Sign and size of $\beta$ are important in the study of constant coefficient system. Particularly, the interesting paper \cite{Sirakov} proved that there exists $0<\bar{\beta_{1}}<\bar{\beta_{2}}$ such that (\ref{question2}) has nontrivial positive solutions for $0<\beta<\bar{\beta_{1}}$ or $\beta>\bar{\beta_{2}}$. For the non-constant potentials case, according to Bohr's corresponding principle, the classical description of the mechanical physics
should be recovered by considering the semiclassical limit, i.e. when the Planck constant goes to
zero. Therefore beside proving existence of solutions to (\ref{question2}), it is physically relevant to study the
asymptotic limit as $\varepsilon\rightarrow0$ of such solutions.\

The semiclassical limit is well understood for the nonlinear Schr\"{o}dinger equation
\begin{align*}
-\varepsilon^{2}\Delta u_{\varepsilon}+Vu_{\varepsilon}=|u_{\varepsilon}|^{p-1}u_{\varepsilon}\ \ \ \mbox{in}\ \mathbb{R}^{N}.
\end{align*}
The semi-classical states of the above equation have different limiting behaviors depend upon $\inf_{\mathbb{R}^{N}}V>0$ or $\inf_{\mathbb{R}^{N}}V=0$, the latter is regarded as the $critical$ $frequency$ case. When $\inf_{\mathbb{R}^{N}}V>0$, P. Rabinowitz \cite{Rab} proved that the above problem has a least energy solution for $\varepsilon>0$ small in the case that $\inf_{\mathbb{R}^{N}}V<\liminf_{|x|\rightarrow\infty}V(x)$. In case $\inf_{\mathbb{R}^{N}}V>0$, solutions concentrating at critical points of the potential $V$ have been construted by topological and variational methods \cite{FW,PF,Rab,Wang} and many others. Motivations can be dated back to Floer and Weinstein in there pioneering work \cite{FW}. However, the semi-classical solutions exhibit quite different characteristic features that the energy governed by the shape of the potential $V$ near the zero set of $V$ in the critical frequency case $\inf_{\mathbb{R}^{N}}V=0$ \cite{W-B,Byeon}.\

The semi-classical case (\ref{question2}) with trapping potentials $a(x)$ and $b(x)$ has been studied in some years. Lin and Wei{\cite{LW}} studied (\ref{question2}) by analyzing least energy non-trivial vector solutions. They studied both of attractive interaction and replusive interaction. Especially, when $\beta>0$, they showed the existence of a least energy. In \cite{EBM}, for small $\beta>0$, Montefusco, Pellacci and Squassina exhibited the existence of nonnegative ground state solutions of (\ref{question2}) concentrating around the local minimum points of the potentials, which are in the same region. Ikoma and Tanaka \cite{NK} also considered the case $\beta>0$. They connected the solutions of (\ref{question2}) with the limiting system with $\lambda_{1}=a(x_{0})$ and $\lambda_{2}=b(x_{0})$ for a fixed $x_{0}\in\mathbb{R}^{N}$. Assume that there exists an open bounded set $\Lambda\subset\mathbb{R}^{N}$ such that
\begin{align*}
\mathop{\inf}\limits_{x_{0}\in\Lambda}m(x_{0})<\mathop{\inf}\limits_{x_{0}\in\partial\Lambda}m(x_{0}),
\end{align*}
where $m(x_{0})$ is the ground state energy level of the limiting system. When $\beta>0$ is relatively small, they constructed a family of solutions to (\ref{question2}) which converge to a positive vector solution. For the large $\beta>0$ case, Shi and Wang \cite{Shi} showed the existence of positive ground state solution of (\ref{question2}) concentrating near the minimum points of potential functions. They also showed that multiple positive concentration solutions exist when the topological structure of the set of minimum points satisfies certain condition. We also refer to \cite{LW,Po,Wei} for the study of (\ref{question2}) when $\beta<0$.
\

Note that in all works mentioned above, they all assumed that $a(x)$ and $b(x)$ are positive bounded away from 0. In this paper, we consider the semi-classical in the critical frequency case
\begin{align*}
\mathop{\inf}\limits_{x\in\mathbb{R}^{N}}a(x)=\mathop{\inf}\limits_{x\in\mathbb{R}^{N}}b(x)=0.
\end{align*}\
\indent For the large $\beta>0$ case, Chen and Zou \cite{CZ} assumed $a(x)$ and $b(x)$ are non-negative and may vanish at someplace and decay to $0$ at infinity. They proved the existence of positive solutions for (\ref{question2}) which concentrate around local minima of the potentials as $\varepsilon\rightarrow0$ under suitable conditions on the behaviour of $a(x)$ and $b(x)$. Lucia and Tang \cite{MZ} obtained synchronized solutions which trapped near zero sets of the potentials for $n-$coupled system under some conditions. In \cite{TZW}, Tang and Wang obtained multi-scale segregated positive solutions which are trapped near the nondegenerate critical points as well as the zero set of potentials by using the Lyapunov-Scjmidt reduction method. Tang and Xie \cite{TX} studied $a(x)$ and $b(x)$ admitted several common isolated connected components by using spectral analysis. They  obtained the existence synchronized positive vector solutions which are trapped in a neighbourhood of the zero set of potentials and also the local maximum points of potentials for $\varepsilon$ small.

\

\noindent {\bf 1.2 Motivation.} In this paper, inspired by Byeon and Wang \cite{W-B}, we are concerned with the asymptotic behavior of positive ground state solutions of (\ref{question2}) in the critical frequency case for $\beta>0$ sufficiently large. Obviously, problem (\ref{question2}) is variational. That is, our main aim is to find nontrivial solutions of problem \eqref{question2} by seeking nontrivial critical points in some weighted Sobolev space of the corresponding energy functional $\mathcal{J}_\varepsilon$ given by
$$
\mathcal{J}_\varepsilon(u,v)=\frac{1}{2}\int_{\mathbb{R}^N}\left[\varepsilon^2(|\nabla u|^2+|\nabla v|^2)+a(x)u^2+b(x)v^2\right]-\frac{1}{4}\int_{\mathbb{R}^{N}}(\mu_{1}v^{4}+2\beta u^{2}v^{2}+\mu_{2}v^{4}).
$$
In the sequel, the potentials $a(x)$ and $b(x)$ may achieve $0$ with homogeneous behaviour or vanish on the closure of an open set, but remain bounded away from $0$ at infinity. Throughout the paper, we assume that

$(f_{1})$ $\mu_{1},\mu_{2},\beta>0$.

$(f_{2})$ $a\in C(\mathbb{R}^{N},[0,+\infty))$ and $0=\inf_{x\in\mathbb{R}^{N}}a(x)<\liminf_{|x|\rightarrow\infty}a(x):=a_\infty<\infty$.

$(f_{3})$ $b\in C(\mathbb{R}^{N},[0,+\infty))$ and $0=\inf_{x\in\mathbb{R}^{N}}b(x)<\liminf_{|x|\rightarrow\infty}b(x):=b_\infty<\infty$.

\noindent One of the difficulties in the study of (\ref{question2}) is that it has semi-trivial solutions of type $(u,0)$ and $(0,v)$, where $u(x)$ and $v(x)$ solve
\begin{align}
&-\varepsilon^{2}\Delta u+a(x)u=\mu_{1}u^{3},\ \ \mbox{in}\  \mathbb{R}^{N}.\label{11}\\
&-\varepsilon^{2}\Delta v+b(x)v=\mu_{2}v^{3}, \ \ \mbox{in}\ \mathbb{R}^{N}.\label{12}
\end{align}
By \cite{W-B}, problem \eqref{11} and \eqref{12} both admit ground state solutions for $\varepsilon$ small enough. Here, we borrow an idea of \cite{Sirakov} to overcome this difficulty. Denote by $U_{1\varepsilon}$ and $V_{1\varepsilon}$ a positive ground state solution of (\ref{11}) and (\ref{12}) respectively.

\

\noindent{\bf 1.3 Main results.}
\noindent Let $\mathcal{A}=\{x\in\mathbb{R}^{N}| a(x)=0\}$, $\mathcal{B}=\{x\in\mathbb{R}^{N}| b(x)=0\}$ and $\Sigma=\mathcal{A}\cup\mathcal{B}$. Throughout this paper, we assume that $0\in\Omega=\mathcal{A}\cap\mathcal{B}\not=\emptyset$.
Let
\begin{align*}
\beta_{0}^\varepsilon=\max\Big\{\mu_{2}\frac{\int_{\mathbb{R}^{N}}V_{1\varepsilon}^{4}}{\int_{\mathbb{R}^{N}}U_{1\varepsilon}^{2}V_{1\varepsilon}^{2}},
\mu_{1}\frac{\int_{\mathbb{R}^{N}}U_{1\varepsilon}^{4}}{\int_{\mathbb{R}^{N}}U_{1\varepsilon}^{2}V_{1\varepsilon}^{2}}\Big\}.\end{align*}

\begin{definition}\cite{W-B}
A function $W : \mathbb{R}^{N}\rightarrow\mathbb{R}$ is called positive $\gamma$-homogeneous if $W(y)>0$ for all $y\neq0$ and $W(ty)=t^{\gamma}W(y)$ for any $t\in[0,+\infty)$ and $y\in\mathbb{R}^{N}$.
\end{definition}
For any given positive $\gamma$-homogeneous function $W$, we consider the following coupled system
\begin{align}\label{limit}
\begin{cases}
-\Delta u+W(x)u=\mu_{1}u^{3}+\beta v^{2}u \ \ \ \ \mbox{in}\ \mathbb{R}^{N},\\
-\Delta v+W(x)v=\mu_{2}v^{3}+\beta u^{2}v \ \ \ \ \ \mbox{in}\ \mathbb{R}^{N},
\end{cases}
\end{align}
which will play as the limit problem of (\ref{question2}) with homogeneous potentials and whose corresponding functional $\mathcal{J}_{W}(u,v)\in C^{1}(H_{W}(\mathbb{R}^{N})\times H_{W}(\mathbb{R}^{N}))$ is defined by
\begin{align*}
\mathcal{J}_{W}(u,v)=\frac{1}{2}\int_{\mathbb{R}^{N}}|\nabla u|^{2}+W(x)u^{2}+|\nabla v|^{2}+W(x)v^{2}-\frac{1}{4}\int_{\mathbb{R}^{N}}\mu_{1}v^{4}+2\beta u^{2}v^{2}+\mu_{2}v^{4}.
\end{align*}
Here $H_{W}$ is the completion of $C^{\infty}_{0}(\mathbb{R}^{N})$ with respect to the norm
\begin{align*}
\|u\|_{H_{W}}=\Big(\int_{\mathbb{R}^{N}}|\nabla u|^{2}+Wu^{2}\Big)^{\frac{1}{2}}.
\end{align*}
Since the potential $W$ is coercive, one can show that, for $\beta>0$,
\begin{align*}
\mathcal{E}_{W}:=\inf\{\mathcal{J}_{W}(u,v)| (u,v)\in H_{W}\times H_{W}\setminus\{(0,0)\},\langle\mathcal{J}_{W}'(u,v),(u,v)\rangle=0\},
\end{align*}
can be attained by a positive ground state solution.

Now, we are in position to state the first result as follows.
\
\begin{theorem}
\label{Theorem 1.1.} {\it Assume that $(f_{1})$-$(f_{3})$ hold and there exists $\gamma>0$ such that for every $x\in\Omega$, one of the following two hypotheses holds
\begin{itemize}
\item [($H_1$)]
$$
\mathop{\lim}\limits_{z\rightarrow x}\frac{a(z)}{|z-x|^{\gamma}}=+\infty,
$$
\item [($H_2$)] for some positive $\gamma$-homogeneous function $W\in C(\mathbb{R}^{N})$, there holds that
$$
\mathop{\lim}\limits_{z\rightarrow x}\frac{a(z)-W(z-x)}{|z-x|^{\gamma}}=0.
$$
\end{itemize}
Additionally, $(H_2)$ holds for at least one point $x\in\Omega$, where there also holds that
\begin{align}\label{bb}
\mathop{\lim}\limits_{z\rightarrow x}\frac{b(z)-W(z-x)}{|z-x|^{\gamma}}=0.
\end{align}
Then for sufficiently small $\varepsilon>0$ and $\beta>\beta_0^\varepsilon$, system (\ref{question2}) has a positive ground state solution $(U_{\varepsilon},V_{\varepsilon})$.\

Moreover, denoting by $\{x_{i}; i\in\mathcal{I}\}$ all points where $(H_2)$ and \eqref{bb} both hold. For $\beta>\beta_0$(see below), there exists $x_{*}\in\{x_{i}; i\in\mathcal{I}\}$ and a positive $\gamma-homogeneous$ function $W_{\ast}\in C(\mathbb{R}^{N})$ such that
\begin{align*}
\mathop{\lim}\limits_{x\rightarrow x_{*}}\frac{a(x)-W_{\ast}(x-x_{*})}{|x-x_{*}|^{\gamma}}=0\  \mbox{and}\ \mathop{\lim}\limits_{x\rightarrow x_{*}}\frac{b(x)-W_{\ast}(x-x_{*})}{|x-x_{*}|^{\gamma}}=0,
\end{align*}
and as $\varepsilon\rightarrow0$,
\begin{align*}
\varepsilon^{-\frac{\gamma}{\gamma+2}}(U_{\varepsilon}(\varepsilon^{\frac{2}{\gamma+2}}\cdot+x_{*}),V_{\varepsilon}(\varepsilon^{\frac{2}{\gamma+2}}\cdot+x_{*}))\rightarrow(\omega,\phi)\ \mbox{ in}\ H_{loc}^{1}(\mathbb{R}^{N})\times H_{loc}^{1}(\mathbb{R}^{N})
\end{align*}
and
\begin{align*}
\mathop{\lim}\limits_{\varepsilon\rightarrow0}\frac{\mathcal{J}_{\varepsilon}(U_{\varepsilon},V_{\varepsilon})}{\varepsilon^{\frac{2}{\gamma+2}(N+2\gamma)}}
=\mathcal{E}_{W_{\ast}}=\mathop{\inf}\limits_{i\in\mathcal{I}}\mathcal{E}_{W_{i}}.
\end{align*}
Here $(w,\phi)$ is a ground state solution of the limit problem \eqref{limit} with $W=W_{x_{*}}$.
}
\end{theorem}
\br\label{remark}

Without loss of generality, assume $0\in\{x_{i};i\in\mathcal{I}\}$. Notice that the conclusion of \cite[Theorem 3.2.3]{W-B} also holds for any least-energy solution. Then it follows from \cite[Theorem 3.2.3]{W-B} that rescaled functions $\varepsilon^{-\frac{\gamma}{\gamma+2}}U_{1\varepsilon}(\e^{\frac{2}{\gamma+2}}\cdot)$ and $\varepsilon^{-\frac{\gamma}{\gamma+2}}V_{1\varepsilon}(\e^{\frac{2}{\gamma+2}}\cdot)$ subconverge to $\omega_{\gamma,\mu_1}$ and $\omega_{\gamma,\mu_2}$ uniformly on $\R^{N}$ respectively as $\e$ goes to zero. Here, for $i=1,2$, $\omega_{\gamma,\mu_i}$ is a least-energy solution of the problem below
$$
-\Delta u+W(x)u=\mu_{i}u^{3}\ \ \ \ \mbox{in}\ \mathbb{R}^{N}.\\
$$
Therefore, $\varepsilon^{-\frac{4\gamma}{\gamma+2}}\int_{\R^{N}}U_{1\e}^{4}(\e^{\frac{2}{\gamma+2}}x)$ and $\varepsilon^{-\frac{4\gamma}{\gamma+2}}\int_{\R^{N}}V_{1\e}^{4}(\e^{\frac{2}{\gamma+2}}x)$ are bounded and
\begin{align*}
&\lim_{\varepsilon\rightarrow0}\varepsilon^{-\frac{4\gamma}{\gamma+2}}\int_{\R^{N}}U_{1\e}^{4}(\e^{\frac{2}{\gamma+2}}x)=\int_{\R^{N}}\omega_{\gamma,\mu_1}^{4}>0,\\
&\lim_{\varepsilon\rightarrow0}\varepsilon^{-\frac{4\gamma}{\gamma+2}}\int_{\R^{N}}V_{1\e}^{4}(\e^{\frac{2}{\gamma+2}}x)=\int_{\R^{N}}\omega_{\gamma,\mu_2}^{4}>0,\\
\end{align*}
and
\begin{align*}
\lim_{\varepsilon\rightarrow0}\int_{\R^{N}}\varepsilon^{-\frac{2\gamma}{\gamma+2}}U_{1\e}^{2}(\e^{\frac{2}{\gamma+2}}x)\cdot\varepsilon^{-\frac{2\gamma}{\gamma+2}}V_{1\e}^2(\e^{\frac{2}{\gamma+2}}x)
=\int_{\R^{N}}\omega_{\gamma,\mu_1}^{2}\omega_{\gamma,\mu_2}^{2}>0.
\end{align*}
 Then
$$
0<\lim_{\varepsilon\rightarrow0}\beta_{0}^\varepsilon:=\beta_{0}=
\max\left\{\frac{\mu_1\int_{\R^{N}}\omega_{\gamma,\mu_1}^{4}}{\int_{\R^{N}}\omega_{\gamma,\mu_1}^{2}\omega_{\gamma,\mu_2}^{2}},
\frac{\mu_2\int_{\R^{N}}\omega_{\gamma,\mu_2}^{4}}{\int_{\R^{N}}\omega_{\gamma,\mu_1}^{2}\omega_{\gamma,\mu_2}^{2}}\right\}<\infty.
$$
\er

\

Next we consider the case that the potential $a$ and $b$ vanish on the nonempty open set with smooth boundary for the nonlinear Schr\"{o}dinger system. That is, we assume that $\mbox{int}(\mathcal{A})$ and $\mbox{int}(\mathcal{B})$ are nonempty with smooth boundary. For the sake of simplicity, $\mbox{int}(\mathcal{A})$ and $\mbox{int}(\mathcal{B})$ are connected.
\br\label{remark2}
{
It follows from \cite[Theorem 3.1]{W-B} that $\varepsilon^{-1}U_{1\varepsilon}$ subconverges pointwise to $U_{0}$ on $\mbox{int}(\mathcal{A})$ and to 0 on $\R^{N}\setminus\mbox{int}(\mathcal{A})$. Similarly, $\varepsilon^{-1}V_{1\varepsilon}$ subconverges pointwise to $V_{0}$ on $\mbox{int}(\mathcal{B})$ and to 0 on $\R^{N}\setminus\mbox{int}(\mathcal{B})$. Here $U_{0}$ is a positive ground state solution of
\begin{align*}
\begin{cases}
-\Delta u=\mu_{1}u^{3}, \ \ \ \  x\in\ \mbox{int}(\mathcal{A}),\\
u\in H_{0}^{1}(\mbox{int}(\mathcal{A}))
\end{cases}
\end{align*}
and
$V_{0}$ is a positive ground state solution of
\begin{align*}
\begin{cases}
-\Delta v=\mu_{2}v^{3}, \ \ \ \  x\in\ \mbox{int}(\mathcal{B}),\\
v\in H_{0}^{1}(\mbox{int}(\mathcal{B})).
\end{cases}
\end{align*} Hence $\e^{-4}\int_{\R^{N}}U_{1\e}^{4}$ and $\e^{-4}\int_{\R^{N}}V_{1\e}^{4}$ are bounded and
\begin{align*}
&\lim_{\e\rightarrow0}\e^{-4}\int_{\R^{N}}U_{1\e}^{4}=\int_{\mathcal{A}}U_{0}^{4}>0,\\
&\lim_{\e\rightarrow0}\e^{-4}\int_{\R^{N}}V_{1\e}^{4}=\int_{\mathcal{B}}V_{0}^{4}>0,
\end{align*}
and
\begin{align*}
\lim_{\e\rightarrow0}\e^{-4}\int_{\R^{N}}U_{1\e}^{2} V_{1\e}^{2}=\int_{\Omega}U_{0}^{2}V_{0}^{2}>0.
\end{align*}
}
Set
\begin{align*}
\beta_{1}=\max\left\{\mu_{2}\frac{\int_{\mathcal{B}}V_{0}^{4}}{\int_{\Omega}U_{0}^{2}V_{0}^{2}},\mu_{1}\frac{\int_{\mathcal{A}}U_{0}^{4}}{\int_{\Omega}U_{0}^{2}V_{0}^{2}}\right\},
\end{align*}
one can get that
$$
\lim_{\e\rightarrow0}\beta_0^\varepsilon=\beta_1.
$$
\er
\noindent Denote by int$(\Omega)$ the interior of subset $\Omega$. Our second result states as following.
\begin{theorem}
\label{Theorem 1.2.} {\it Let $\beta>\beta_1$. Assume that $(f_{1})$-$(f_{3})$ hold and int$(\Omega)$ is a nonempty open set with smooth boundary,  then, for $\varepsilon>0$ small, system (\ref{question2}) has a positive ground state solution $(U_{\varepsilon},V_{\varepsilon})$. Moreover, if $\mu_1,\mu_2>0$ small, there exists a ground state solution $(w,z)\in H^{1}(\mathbb{R}^{N})\times H^{1}(\mathbb{R}^{N})$ of the following problem
\begin{align*}
\begin{cases}
-\Delta w=\mu_{1}w^{3}+\beta z^{2}w, \   x\in\mbox{int}(\Omega),\\
-\Delta z=\mu_{2}z^{3}+\beta w^{2}z, \  x\in\mbox{int}(\Omega),\\
w, z\in H_{0}^{1}(\mbox{int}(\Omega))
\end{cases}
\end{align*}
such that, as $\varepsilon\rightarrow0$,
\begin{align*}
\varepsilon^{-1}(U_{\varepsilon},V_{\varepsilon})\rightarrow(w,z)\ \ \mbox{in}\ H^{1}(\mathbb{R}^{N})\times H^{1}(\mathbb{R}^{N}).
\end{align*}
}
\end{theorem}

For the proof of our theorems, we shall consider an equivalent system to (\ref{question2}). By the scaling $a(x)=a(\varepsilon x)$, $b(x)=b(\varepsilon x)$, we then deduce that\
\begin{align*}
\label{question3}
\begin{cases}
-\Delta u+a(\varepsilon x)u=\mu_{1}u^{3}+\beta uv^{2} \ \ \ \ x\in\mathbb{R}^{N},\\
-\Delta v+b(\varepsilon x)v=\mu_{1}v^{3}+\beta u^{2}v \ \ \ \ x\in\mathbb{R}^{N}.\tag{$\mathcal{P}_{\varepsilon}$}
\end{cases}
\end{align*}
Thus, in the sequel we focus only on system (\ref{question3}).

The  rest of this paper is organized as follows. In Section 2, we introduce a variational setting of our problem and present some preliminary results about the extended problem. Section 3 is devoted to the study of positive ground state solutions of problem (\ref{question3}). Section 4 is devoted to the proof of Theorem \ref{Theorem 1.1.}. Section 5 is devoted to the proof of Theorem \ref{Theorem 1.2.}.

\s {Preliminaries}
\renewcommand{\theequation}{2.\arabic{equation}}
To prove the main results, we use the following notations.\

\

$\bullet$ $\mathbb{H} :=H^{1}(\mathbb{R}^{N})\times H^{1}(\mathbb{R}^{N})$ with norm
\begin{align*}
\|(u,v)\|_{\mathbb{H}}=(\|u\|_{H^{1}}^{2}+\|v\|_{H^{1}}^{2})^{1/2};
\end{align*}

$\bullet$ $\mathbb{H}_{loc} :=H_{loc}^{1}(\mathbb{R}^{N})\times H_{loc}^{1}(\mathbb{R}^{N})$

$\bullet$ $\mathbb{L}^{q}(\mathbb{R}^{N}) :=L^{q}(\mathbb{R}^{N})\times L^{q}(\mathbb{R}^{N})$;

$\bullet$ $\mathbb{L}^{q}_{loc}(\mathbb{R}^{N}) :=L^{q}_{loc}(\mathbb{R}^{N})\times L^{q}_{loc}(\mathbb{R}^{N})$;

$\bullet$ $\mathbb{C}^{\infty}_{0}(\mathbb{R}^{N}) :=C^{\infty}_{0}(\mathbb{R}^{N})\times C^{\infty}_{0}(\mathbb{R}^{N})$;

$\bullet$ $2^{*}=6$ when $N=3$, and $2^{*}=+\infty$ when $N=1,2$.

\
For $\varepsilon>0$, we define the rescaled weighted Sobolev space $H_{a,\varepsilon}$ by
\begin{align*}
H_{a,\varepsilon}=\left\{u\in H^{1}(\mathbb{R}^{N})|\int_{\mathbb{R}^{N}}a(\varepsilon x)u^{2}<+\infty\right\},
\end{align*}
endowed with the norm
\begin{align*}
\|u\|_{a,\varepsilon}=\Big(\int_{\mathbb{R}^{N}}|\nabla u|^{2}+a(\varepsilon x)u^{2}\Big)^{\frac{1}{2}}.
\end{align*}
Define $E :=H_{a,\varepsilon}\times H_{b,\varepsilon}$ with norm
\begin{align*}
\|(u,v)\|_{E}=(\|u\|_{a,\varepsilon}^{2}+\|v\|_{b,\varepsilon}^{2})^{1/2}.
\end{align*}
For any $(u,v)\in E$, set
\begin{align*}
F(u,v)=\int_{\mathbb{R}^{N}}\mu_{1}u^{4}+2\beta u^{2}v^{2}+\mu_{2}v^{4}
\end{align*}
and the energy functional for (\ref{question3}) is defined by
\begin{align*}
I_{\varepsilon}(u,v) :=\frac{1}{2}\int_{\mathbb{R}^{N}}|\nabla u|^{2}+a(\varepsilon x)u^{2}+|\nabla v|^{2}+b(\varepsilon x)v^{2}-\frac{1}{4}F(u,v).
\end{align*}

\
We first illustrate the space $E$ can be embedded continuously into the Sobolev space $\mathbb{H}$ for fixed $\varepsilon>0$, even though the potential $a$ and $b$ have a nontrivial set of zeros. Similar to \cite{Schaftingen}, we have

\
\begin{lemma}
\label{Lemma 2.1.} {\it Let $a :\mathbb{R}^{N}\rightarrow[0,+\infty)$. If $\liminf_{|x|\rightarrow\infty}a(x)>0$ then for every $\varepsilon>0$, there exists a constant $C_{\varepsilon}>0$ such that for $u\in H_{a,\varepsilon}(\mathbb{R}^{N})$, $u\in H^{1}(\mathbb{R}^{N})$ and
\begin{align*}
\int_{\mathbb{R}^{N}}|\nabla u|^{2}+|u|^{2}\leq C_{\varepsilon}\int_{\mathbb{R}^{N}}|\nabla u|^{2}+a(\varepsilon x)u^{2}.
\end{align*}
}
\end{lemma}

\

Similar result holds for $H_{b,\varepsilon}$. Lemma \ref{Lemma 2.1.} implies that the space $E$ embeds into $\mathbb{H}$ continuously. Moreover the space $E$ can be continuously embedded in $\mathbb{L}^{q}(\mathbb{R}^{N})$ when $q\in[2,2^{*}]$. It then follows that the functional $I_{\varepsilon}$
is well defined in $E$ and $I_{\varepsilon}\in C^{1}(E,\mathbb{R})$.\

\bd
$(u,v)$ is called a ground state solution of problem (\ref{question3}), if $I_{\varepsilon}(u,v)$ is the least among all nontrival critical points of $I_{\varepsilon}$. Namely, $(u,v)$ has the least energy among nontrivial solutions.
\ed
 A natural method to search the groundstate is to minimize the functional $I_{\varepsilon}$ on the Nehari manifold of problem (\ref{question3}) which is defined by
\begin{align*}
\mathcal{N}_{\varepsilon}:=\left\{(u,v)\in E\setminus\left\{(0,0)\right\}|\langle I_{\varepsilon}'(u,v),(u,v)\rangle=0\right\}.
\end{align*}
The corresponding groundstate energy is described as
\begin{align*}
c_{\varepsilon}:=\mathop{\inf}\limits_{(u,v)\in\mathcal{N}_{\varepsilon}}I_{\varepsilon}(u,v).
\end{align*}\

\begin{lemma}
\label{Lemma 2.2.} {\it For given $\varepsilon>0$, the ground state energy $c_{\varepsilon}$ is positive and $\mathcal{N}_{\varepsilon}$ is a manifold of class of $C^{1}$. Moreover, if $(u,v)\in\mathcal{N}_{\varepsilon}$ is a critical point of $I_{\varepsilon}|_{\mathcal{N}_{\varepsilon}}$, then $I_{\varepsilon}'(u,v)=0$ in $E^\ast$.}
\end{lemma}

\noindent{\it Proof.} For any fixed $\varepsilon>0$ and $(u,v)\in E\setminus\left\{(0,0)\right\}$, we have that
\begin{align*}
t(u,v)\in\mathcal{N}_{\varepsilon}\Leftrightarrow t^{2}\|(u,v)\|_{E}^{2}=t^{4}F(u,v).
\end{align*}
As a consequence, for all $(u,v)\in E\setminus\left\{(0,0)\right\}$, there exists a unique $t>0$ such that $t(u,v)\in\mathcal{N}_{\varepsilon}$. Moreover, since $F(u,v)$ is homogeneous with degree 4, there is $\rho>0$ such that
\begin{align*}
\|(u,v)\|_{E}^{2}\geq\rho>0,\ \ \ \forall(u,v)\in\mathcal{N}_{\varepsilon}.
\end{align*}
Hence, for any $(u,v)\in\mathcal{N}_{\varepsilon}$, we have
\begin{align*}
I_{\varepsilon}(u,v)=\frac{1}{4}\|(u,v)\|_{H_{\varepsilon}}^{2}\geq\frac{1}{4}\rho>0,
\end{align*}
thus, $c_{\varepsilon}>0$.
We define $\mathfrak{g}_{\varepsilon}(u,v)=\langle I_{\varepsilon}'(u,v),(u,v)\rangle$ for each $(u,v)\in\mathcal{N}_{\varepsilon}$,
\begin{align}\label{3}
\langle\mathfrak{g}_{\varepsilon}'(u,v),(u,v)\rangle=-2\|(u,v)\|_{E}^{2}\leq-2\rho<0,
\end{align}
which implies that $\mathcal{N}_{\varepsilon}$ is a smooth complete manifold of codimension one in $E$.\
Moreover, if $(u,v)\in\mathcal{N}_{\varepsilon}$ is a critical point of the restricted functional $I_{\varepsilon}|_{\mathcal{N}_{\varepsilon}}$, then there exists $\lambda_{\varepsilon}\in\mathbb{R}$, such that
\begin{align}\label{4}
I_{\varepsilon}'(u,v)=\lambda_{\varepsilon}\mathfrak{g}_{\varepsilon}'(u,v)\,\,\mbox{in}\,\,E^\ast.
\end{align}
Then, we have
\begin{align*}
\lambda_{\varepsilon}\langle\mathfrak{g}_{\varepsilon}'(u,v),(u,v)\rangle=\langle I_{\varepsilon}'(u,v),(u,v)\rangle=0
\end{align*}
and deduce by (\ref{3}) that $\lambda_{\varepsilon}=0$. The conclusions follows then from (\ref{4}).
\qed

\s{Existence of positive ground state solutions of (\ref{question3})}
\renewcommand{\theequation}{3.\arabic{equation}}
\begin{lemma}
\label{Lemma 3.1.}{\it(Mountain-Pass geometry)} {\it Assume that ($f_1$) holds, then the functional $I_{\varepsilon}$ satisfies the following conditions.\\
$(i)$ There exists a positive constant $r>0$ such that $I_{\varepsilon}(u,v)>0$ for $\|(u,v)\|_E=r$;\\
$(ii)$ There exists $(e_{1},e_{2})\in E$ with $\|(e_{1},e_{2})\|_{E}>r$ such that $I_{\varepsilon}(e_{1},e_{2})<0$.\
}
\end{lemma}

\

\noindent{\it Proof.}
\begin{align*}
I_{\varepsilon}(u,v)&=\frac{1}{2}\|(u,v)\|_{E}^{2}-F(u,v);\\
&\geq\frac{1}{2}\|(u,v)\|_{E}^{2}-C\|(u,v)\|_{E}^{4}.
\end{align*}
Hence, there exists $r>0$ such that
\begin{align*}
\mathop{\inf}\limits_{\|(u,v)\|_{E}=r}I_{\varepsilon}(u,v)>0.
\end{align*}
For any $(u,v)\in E\setminus{(0,0)}$ and $t>0$,
\begin{align*}
I_{\varepsilon}(tu,tv)=\frac{t^{2}}{2}\|(u,v)\|_{E}^{2}-\frac{t^{4}}{4}F(u,v),
\end{align*}
which implies that $I_{\varepsilon}(tu,tv)\rightarrow-\infty$ as $t\rightarrow+\infty$.\  $\Box$

 That is, $I_{\varepsilon}$ has a mountaion-pass structure. From Lemma \ref{Lemma 3.1.}, one can apply the Mountain-Pass Theorem{\cite{AP}} without $(PS)_{c}$ condition, and it follows that for any small $\varepsilon>0$, there exists a $(PS)_{c}$ sequence $\left\{(u_{n},v_{n})\right\}\subset E$ (with $c=c_{\varepsilon}^{*}$ defined below) such that
\begin{align*}
&I_{\varepsilon}(u_{n},v_{n})\rightarrow c_{\varepsilon}^{*}=\mathop{\inf}\limits_{\gamma\in\Gamma_{\varepsilon}}\mathop{\max}\limits_{0\leq t\leq1}I_{\varepsilon}(\gamma(t)),\\
&I_{\varepsilon}'(u_{n},v_{n})\rightarrow0, \ \ \mbox{as}\ n\rightarrow\infty,
\end{align*}
where $\Gamma_{\varepsilon}=\left\{\gamma\in C([0,1],E):I_{\varepsilon}(\gamma(0))=0, I_{\varepsilon}(\gamma(1))<0\right\}$.
\qed
\
\begin{lemma}
\label{Lemma 3.2.} {\it $c_{\varepsilon}=c_{\varepsilon}^{*}$.}
\end{lemma}

\noindent{\it Proof.} Define
\begin{align*}
c_{\varepsilon}^{**}&:=\mathop{\inf}\limits_{(u,v)\in E\setminus\left\{(0,0)\right\}}\mathop{\max}\limits_{t\geq0}I_{\varepsilon}(tu,tv).
\end{align*}
For any $(u,v)\in E\setminus\{(0,0)\}$, there exists a unique $t_{u,v}>0$ such that $(t_{u,v}u,t_{u,v}v)\in\mathcal{N}_{\varepsilon}$, and
\begin{align*}
I_{\varepsilon}(t_{u,v}u,t_{u,v}v)=\mathop{\max}\limits_{t\geq0}I_{\varepsilon}(tu,tv),
\end{align*}
which follows that $c_{\varepsilon}=c_{\varepsilon}^{**}$. There exists $t_{0}>0$ large enough such that $I_{\varepsilon}(t_{0}u,t_{0}v)<0$, and
\begin{align*}
c_{\varepsilon}^{*}\leq\mathop{\max}\limits_{t\geq0}I_{\varepsilon}(tt_{0}u,tt_{0}v),
\end{align*}
$c_{\varepsilon}^{*}\leq c_{\varepsilon}^{**}=c_{\varepsilon}$ can be obtained. To show $c_{\varepsilon}^{*}\geq c_{\varepsilon}$, we only need to prove that for any $\gamma\in\Gamma_{\varepsilon}$, $\gamma([0,1])\cap\mathcal{N}_{\varepsilon}\neq\emptyset$.  For any $(u,v)\in E$ with $0<\|(u,v)\|_{E}\ll 1$,
\begin{align*}
\int|\nabla u|^{2}+a(\varepsilon x)u^{2}+|\nabla v|^{2}+b(\varepsilon x)v^{2}>\int\mu_{1}u^{4}+2\beta u^{2}v^{2}+\mu_{2}v^{4}.
\end{align*}
Since $I_{\varepsilon}(\gamma(0))=0$, $\gamma$ is continuous in $t\in[0,1]$, $I_{\varepsilon}(\gamma(t))>0$ for $0<t<\delta$. Thus, if set $\gamma_{t}=(u_{t},v_{t})$,
\begin{align*}
\int|\nabla u_{t}|^{2}+a(\varepsilon x)u_{t}^{2}+|\nabla v_{t}|^{2}+b(\varepsilon x)v_{t}^{2}>\int\mu_{1}u_{t}^{4}+2\beta u_{t}^{2}v_{t}^{2}+\mu_{2}v_{t}^{4}
\end{align*}
holds for sufficiently small $t$. If $\gamma([0,1])\cap\mathcal{N}_{\varepsilon}=\emptyset$, then
\begin{align*}
\int|\nabla u_{1}|^{2}+a(\varepsilon x)u_{1}^{2}+|\nabla v_{1}|^{2}+b(\varepsilon x)v_{1}^{2}>\int\mu_{1}u_{1}^{4}+2\beta u_{1}^{2}v_{1}^{2}+\mu_{2}v_{1}^{4}.
\end{align*}
This is a contradiction with $I_{\varepsilon}(\gamma(1))<0$.
\qed

\

In order to obtain the existence of ground state solutions for (\ref{question3}), we prove some compactness lemma for the functional $I_{\varepsilon}$ to analysis the Palais-Smale sequence properties for the functional $I_{\varepsilon}$.

\

Now, we consider the limit problem
\begin{align*}
\begin{cases}
-\Delta u+\tau u=\mu_{1}u^{3}+\beta v^{2}u \ \ \ \ \mbox{in}\ \mathbb{R}^{N},\\
-\Delta v+\theta v=\mu_{2}v^{3}+\beta u^{2}v \ \ \ \ \ \mbox{in}\ \mathbb{R}^{N},
\end{cases}
\end{align*}
where $(\tau,\theta)\in(0,+\infty)\times(0,\infty)$. The corresponding energy functional and Nehari manifold are defined by
\begin{align*}
I_{\tau\theta}(u,v)=\frac{1}{2}\int_{\mathbb{R}^{N}}|\nabla u|^{2}+\tau u^{2}+|\nabla v|^{2}+\theta v^{2}-\frac{1}{4}F(u,v)
\end{align*}
and
\begin{align*}
\mathcal{N}_{\tau\theta}=\{(u,v)\in \mathbb{H}\setminus\{(0,0)\}|\langle I_{\tau\theta}'(u,v),(u,v)\rangle=0\}.
\end{align*}
Define
\begin{align*}
c_{\tau\theta}=\mathop{\inf}\limits_{(u,v)\in\mathcal{N}_{\tau\theta}}I_{\tau\theta}(u,v).
\end{align*}
One can show that $c_{\tau\theta}$ can be achieved by some $(u,v)\in \mathbb{H}\setminus\{(0,0)\}$ for any $\tau,\theta,\beta>0$.
\begin{lemma}
\label{Lemma 3.3 }{\it The map $c:(0,+\infty)\times(0,+\infty)\rightarrow\mathbb{R}$; $(\tau,\theta)\mapsto c_{\tau\theta}$, is non-decreasing in $\tau$ and $\theta$.}
\end{lemma}

\noindent{\it Proof.} The proof is similar to the proofs of Lemma 2.3 in \cite{Shi} and we omit it here.
\qed
\
\begin{lemma}
\label{Lemma 3.4.} {\it Assume  that $(f_2)$-$(f_3)$ hold and $\Omega\neq\varnothing$. For any fixed $\varepsilon>0$, $I_{\varepsilon}$ satisfies Palais-Smale condition at level $c<c_{a_{\infty}b_{\infty}}$.}
\end{lemma}

\noindent{\it Proof.} 
Let  $c<c_{a_{\infty}b_{\infty}}$, $\{(u_{n},v_{n})\}\subset E$ be a sequence such that $I_{\varepsilon}(u_{n},v_{n})\rightarrow c$, and $I_{\varepsilon}'(u_{n},v_{n})\rightarrow 0$. A standard argument leads us to the fact that ${(u_{n},v_{n})}$ is bounded in $E$, which implies that
\begin{align*}
(u_{n},v_{n})\rightharpoonup(u,v)\ \mbox{weakly in}\ E,\ I_{\varepsilon}'(u,v)=0,\ I_{\varepsilon}(u,v)\geq0.
\end{align*}
Set $w_{n}=u_{n}-u,z_{n}=v_{n}-v$, and $(w_{n},z_{n})\rightharpoonup0$. By Brezis-Lieb Lemma {\cite{Willem}}, it follows that,
\begin{align*}
&I_{\varepsilon}(w_{n},z_{n})=I_{\varepsilon}(u_{n},v_{n})-I_{\varepsilon}(u,v)+o(1)\leq c+o(1),\\
&I_{\varepsilon}'(w_{n},z_{n})\rightarrow0\ \ \mbox{as}\ n\rightarrow\infty.
\end{align*}
There exists $t_{n}>0$ such that $t_{n}(w_{n},v_{n})\in \mathcal{N}_{a_{\infty},b_{\infty}}$, from which it follows that
\begin{align}\label{5}
t_{n}^{2}\int_{\mathbb{R}^{N}}|\nabla w_{n}|^{2}+a_{\infty}w_{n}^{2}+|\nabla z_{n}|^{2}+b_{\infty}z_{n}^{2}=t_{n}^{4}\int_{\mathbb{R}^{N}}\mu_{1}w_{n}^{4}+2\beta w_{n}^{2}z_{n}^{2}+\mu_{2}z_{n}^{4},
\end{align}
Suppose that $(w_{n},z_{n})\nrightarrow(0,0)$ in $E$, we first claim that
\begin{align}\label{6}
\limsup_{n\rightarrow\infty}t_{n}\leq1.
\end{align}
Otherwise, we can assume that $t_{n}\geq1+\delta$ for $\delta>0$. From $I_{\varepsilon}'(w_{n},z_{n})=o(1)$, we have that
\begin{align}\label{7}
\int_{\mathbb{R}^{N}}|\nabla w_{n}|^{2}+a(\varepsilon x)w_{n}^{2}+|\nabla z_{n}|^{2}+b(\varepsilon x)z_{n}^{2}=\int_{\mathbb{R}^{N}}\mu_{1}w_{n}^{4}+2\beta w_{n}^{2}z_{n}^{2}+\mu_{2}z_{n}^{4}+o(1).
\end{align}
According to $(\ref{5})$ and $(\ref{7})$ , we have
\begin{align}\label{8}
\int_{\mathbb{R}^{N}}(a_{\infty}-a(\varepsilon x))w_{n}^{2}+(b_{\infty}-b(\varepsilon x))z_{n}^{2}=(t_{n}^{2}-1)\int_{\mathbb{R}^{N}}\mu_{1}w_{n}^{4}+2\beta w_{n}^{2}z_{n}^{2}+\mu_{2}z_{n}^{4}+o(1)
\end{align}
For any $\varepsilon>0$, there exists $R=R(\varepsilon)$ such that
\begin{align}\label{9}
a(\varepsilon x)>a_{\infty}-\varepsilon\ \  \mbox{and}\ \ b(\varepsilon x)>b_{\infty}-\varepsilon\ \ \mbox{for}\  |x|\geq R.
\end{align}
Since $(w_{n},z_{n})\rightarrow0$ in $\mathbb{L}^{2}_{B_{R}(0)}$ and $(\ref{9})$, there exists $C>0$ such that
\begin{align*}
(t_{n}^{2}-1)\int_{\mathbb{R}^{N}}\mu_{1}w_{n}^{4}+2\beta w_{n}^{2}z_{n}^{2}+\mu_{2}z_{n}^{4}\leq C\varepsilon+o(1)
\end{align*}
Since $(w_{n},z_{n})\nrightarrow(0,0)$, there exists a sequence $\{y_{n}\}\in\mathbb{R}^{N}$ and positive constants $r,\zeta>0$,
\begin{align}\label{10}
\mathop{\liminf}\limits_{n\rightarrow\infty}\int_{B_{r}(y_{n})}(w_{n}^{2}+z_{n}^{2})\geq\zeta.
\end{align}
Set $(\tilde{w}_{n}(x),\tilde{z}_{n}(x))=(u(x+y_{n}),v(x+y_{n}))$. Then $\left\{(\tilde{w}_{n},\tilde{z}_{n})\right\}$ is bounded and $(\tilde{w}_{n},\tilde{z}_{n})\rightharpoonup(u,v)$ in $E$, $(\tilde{w}_{n},\tilde{z}_{n})\rightarrow(u,v)$ in $\mathbb{L}_{loc}^{2}(\mathbb{R}^{N})$. Moreover, there exists a subset $\Lambda\subset\mathbb{R}^{N}$ such that $(u,v)\neq(0,0)$ $a.e.$ in $\Lambda$, then it follows from Fatou's lemma that
\begin{align*}
0<[(1+\delta)^{2}-1]\int_{\Lambda}\mu_{1}u^{4}+2\beta u^{2}v^{2}+\mu_{2}v^{4}\leq C\varepsilon,
\end{align*}
for any $\varepsilon>0$, which yields a contradiction. Thus (\ref{6}) holds. Using the fact that $t_{n}(w_{n},v_{n})\in \mathcal{N}_{a_{\infty},b_{\infty}}$, and (\ref{6})
\begin{align*}
c_{a_{\infty}b_{\infty}}&\leq I_{a_{\infty}b_{\infty}}(t_{n}w_{n},t_{n}z_{n})-\frac{1}{4}\langle I_{\infty}'(t_{n}(w_{n},v_{n})),t_{n}(w_{n},v_{n})\rangle\\
&=\frac{t_{n}^{2}}{4}\int_{\mathbb{R}^{N}}|\nabla w_{n}|^{2}+a_{\infty}w_{n}^{2}+|\nabla z_{n}|^{2}+b_{\infty}z_{n}^{2}\\
&\leq\frac{1}{4}\int_{\mathbb{R}^{N}}|\nabla w_{n}|^{2}+a(\varepsilon x)w_{n}^{2}+|\nabla z_{n}|^{2}+b(\varepsilon x)z_{n}^{2}+o(1)\\
&=I_{\varepsilon}(w_{n},z_{n})+o(1)\leq c+o(1)
\end{align*}
Let $n\rightarrow\infty$ in the above inequality, we have $c_{a_{\infty}b_{\infty}}\leq c$ which is impossible. Hence $\limsup_{n\rightarrow\infty}t_{n}\leq1$ cannot happen. It implies $(u_{n},v_{n})\rightarrow(u,v)$ in $E$.
\qed
\
\begin{lemma}
\label{Lemma 3.5.} {\it $c_{\varepsilon}<c_{a_{\infty}b_{\infty}}$} for sufficiently small $\varepsilon>0$.
\end{lemma}

\noindent{\it Proof.} Taking any $\tau\in(0,a_\infty), \theta\in(0,b_\infty)$, let $(u,v)$ be a ground state solution for $I_{\tau\theta}$. Choose $\eta\in C_{0}^{\infty}(\mathbb{R}^{N},[0,1])$ such that $\eta(x)=1$, for $|x|\leq\frac{1}{2}$, and $\eta(x)=0$, for $|x|\geq1$.  For any $n\in\mathbb{N}$, define $(u_{n},v_{n})=\eta(\frac{\cdot}{n})(u,v)$ and $(w_{n},z_{n})=t_{n}(u_{n},v_{n})$ such that $(w_{n},z_{n})\in\mathcal{N}_{\tau\theta}$. Since $(u_{n},v_{n})\rightarrow(u,v)$ in $\mathbb{H}$ as $n\rightarrow\infty$, and $(u,v)\in\mathcal{N}_{\tau\theta}$, it is easy to see that $\lim_{n\rightarrow\infty}t_{n}=1$. Consequently, as $n\rightarrow\infty$, we have
\begin{align*}
I_{\tau\theta}(w_{n},z_{n})=I_{\tau\theta}(t_{n}(u_{n},v_{n}))\rightarrow c_{\tau\theta}< c_{a_{\infty}b_{\infty}}.
\end{align*}
Now, we have $(w_{n},z_{n})\rightarrow(u,v)$ in $\mathbb{H}$ and $(w_{n},z_{n})\in\mathcal{N}_{\tau\theta}$, which implies that
\begin{align*}
I_{\tau\theta}(w_{n},z_{n})=\mathop{\max}\limits_{t\geq0}I_{\tau\theta}(tw_{n},tz_{n}).
\end{align*}
By assumption $(f_2)$-$(f_3)$ and $0\in\Omega=\mathcal{A}\cap\mathcal{B}$, for any $\tau,\theta>0$, there exists $R>0$ small such that $a(x)<\tau,b(x)<\theta$ for $|x|<R$ and
\begin{align*}
a(\varepsilon x)<\tau,\ \ \ b(\varepsilon x)<\theta,\ \ \ \mbox{for}\ |x|<\frac{R}{\varepsilon}.
\end{align*}
Define
\begin{align*}
\mathbb{H}_{R/\varepsilon} :=H_{0}^{1}(B_{R/\varepsilon}(0))\times H_{0}^{1}(B_{R/\varepsilon}(0)).\\
c_{\tau\theta}^{R/\varepsilon} :=\mathop{\inf}\limits_{(u,v)\in\mathbb{H}_{R/\varepsilon}\setminus\{(0,0)\}}\mathop{\max}\limits_{t\geq0}I_{\tau\theta}(tu,tv),
\end{align*}
then
\begin{align*}
c_{\varepsilon}&=\mathop{\inf}\limits_{(u,v)\in E\setminus\{(0,0)\}}\mathop{\max}\limits_{t\geq0}I_{\varepsilon}(tu,tv)\leq\mathop{\inf}\limits_{(u,v)\in\mathbb{H}_{R/\varepsilon}\setminus\{(0,0)\}}\mathop{\max}\limits_{t\geq0}I_{\varepsilon}(tu,tv)\\
&\leq\mathop{\inf}\limits_{(u,v)\in\mathbb{H}_{R/\varepsilon}\setminus\{(0,0)\}}\mathop{\max}\limits_{t\geq0}I_{\tau\theta}(tu,tv)=c_{\tau\theta}^{R/\varepsilon}.
\end{align*}
Then for any fixed $n$ with $I_{\tau\theta}(w_{n},z_{n})<c_{a_{\infty}b_{\infty}}$ and sufficiently small $\varepsilon>0$ such that $R/\varepsilon>n$, we have $(w_{n},z_{n})\in \mathbb{H}_{R/\varepsilon}$ and
\begin{align*}
c_{\tau\theta}^{R/\varepsilon}=\mathop{\inf}\limits_{(u,v)\in\mathbb{H}_{R/\varepsilon}\setminus\{(0,0)\}}\mathop{\max}\limits_{t\geq0}I_{\tau\theta}(tu,tv)
\leq\mathop{\max}\limits_{t\geq0}I_{\tau\theta}(tw_{n},tz_{n})=I_{\tau\theta}(w_{n},z_{n})<c_{a_{\infty}b_{\infty}}.
\end{align*}
Then the conclusion follows $c_{\varepsilon}<c_{a_{\infty}b_{\infty}}$ for sufficiently small $\varepsilon>0$.
\qed
\
\begin{lemma}
\label{Lemma 3.6.}{\it Assume $(f_1)$ holds, then for $\beta>\beta_0^\varepsilon$ and sufficiently small $\varepsilon>0$, (\ref{question3}) has a nontrivial ground state solution.}
\end{lemma}

\noindent{\it Proof.} From Lemmas \ref{Lemma 3.4.}-\ref{Lemma 3.6.}, one sees that problem (\ref{question3}) has a ground state solution $(u_{\e},v_{\e})$ for sufficiently small $\e$. Thus, $(U_{\e},V_{\e})=(u_{\e}(\frac{\cdot}{\e}),v_{\e}(\frac{\cdot}{\e}))$ is a ground state solution of problem (\ref{question2}) for sufficiently small $\e$. It is sufficient to prove that for sufficiently small $\varepsilon>0$,
\begin{align}\label{13}
c_{\varepsilon}<\min\{I_{\varepsilon}(U_{1\varepsilon}(\e\cdot),0),I_{\varepsilon}(0,V_{1\varepsilon}(\e\cdot)\}.
\end{align}
First, as in Lemma 3.3 of \cite{Sirakov}, we know that for $\beta>0$,
\begin{align*}
c_{\varepsilon}=\mathop{\inf}\limits_{(u,v)\in E\setminus\{(0,0)\}}\mathcal{K}_{\varepsilon}(u,v),
\end{align*}
where
\begin{align*}
\mathcal{K}_{\varepsilon}(u,v)=\frac{(\|u\|_{a,\varepsilon}^{2}+\|v\|_{b,\varepsilon}^{2})^{2}}{4F(u,v)}.
\end{align*}
For $(s,t)\in\Gamma=\{(s,t): s\geq0,t\geq0,(s,t)\neq(0,0)\}$, we define a function
\begin{align*}
f_{\e}(s,t)=\mathcal{K}_{\varepsilon}(\sqrt{s}U_{1\varepsilon}(\e\cdot),\sqrt{t}V_{1\varepsilon}(\e\cdot))=\frac{(s\|U_{1\varepsilon}(\e\cdot)\|_{a,\varepsilon}^{2}+t\|V_{1\varepsilon}(\e\cdot)\|_{b,\varepsilon}^{2})^{2}}{4\int_{\mathbb{R}^{N}}
\mu_{1}s^{2}U_{1\varepsilon}^{4}(\e x)+2\beta stU_{1\varepsilon}^{2}(\e x){V_{1\varepsilon}^{2}(\e x)}+\mu_{2}t^{2}V_{1\varepsilon}^{4}(\e x)}.
\end{align*}
It follows that
\begin{align*}
f_{\e}(s,0)=\frac{\|U_{1\varepsilon}(\e\cdot)\|_{a,\varepsilon}^{4}}{4\int_{\mathbb{R}^{N}}\mu_{1}U_{1\varepsilon}^{4}(\e x)}=I_{\varepsilon}(U_{1\varepsilon}(\e\cdot),0),\\
f_{\e}(0,t)=\frac{\|V_{1\varepsilon}(\e\cdot)\|_{b,\varepsilon}^{4}}{4\int_{\mathbb{R}^{N}}\mu_{2}V_{1\varepsilon}^{4}(\e x)}=I_{\varepsilon}(0,V_{1\varepsilon}(\e\cdot)).
\end{align*}
Thus, to guarantee (\ref{13}) it is sufficient to prove that $f_\varepsilon$ does not attain its minimum over $\Gamma$ on the line $s=0$ or $t=0$. The function $f_\varepsilon$ is a fraction of two quadratic forms in $(s,t)$, and elementary analysis shows that the quantity
\begin{align*}
\frac{(as+bt)^{2}}{cs^{2}+2dst+et^{2}},\ \ \ a,b,c,d,e>0
\end{align*}
does not attain its minimum in $\Gamma$ on the axes if and only if
\begin{align}\label{14}
ad-bc>0,\ \ \ \ \ bd-ae>0
\end{align}
and then the minimum is attained at $(ad-bc,bd-ae)$. Clearly, to verify (\ref{13}), it suffices to show
\begin{align*}
\beta>\mu_{2}\frac{\int_{\mathbb{R}^{N}}V_{1\varepsilon}^{4}}{\int_{\mathbb{R}^{N}}U_{1\varepsilon}^{2}V_{1\varepsilon}^{2}}\ \ \mbox{and}\ \ \beta>\mu_{1}\frac{\int_{\mathbb{R}^{N}}U_{1\varepsilon}^{4}}{\int_{\mathbb{R}^{N}}U_{1\varepsilon}^{2}V_{1\varepsilon}^{2}}.
\end{align*}
Thus, if $\beta>\beta_{0}^\varepsilon$, then (\ref{13}) holds sufficiently small $\varepsilon>0$.
\qed
\
\begin{proposition}
\label{Proposition 3.7.} {\it Assume that $(f_1)$-$(f_3)$ hold, then for $\beta>\beta_{0}^\varepsilon$ and small $\varepsilon>0$, system (\ref{question3}) has a positive ground state solution.}
\end{proposition}

\noindent{\it Proof.} From Lemma \ref{Lemma 3.1.} the functional $I_{\varepsilon}$ satisfies a mountain pass geometry. By Lemma \ref{Lemma 3.2.} and Mountain Pass Theorem \cite{AP}, there exists a $(PS)_{c_{\varepsilon}}$ sequence $\{(u_{n},v_{n})\}\subset E$ such that $I_{\varepsilon}(u_{n},v_{n})\rightarrow c_{\varepsilon}$ and $I_{\varepsilon}'(u_{n},v_{n})\rightarrow0$. Then standard argument shows that $\{(u_{n},v_{n})\}$ is bounded. Up to a subsequence, Lemma \ref{Lemma 3.5.} implies that there exists $(u_{\varepsilon},v_{\varepsilon})$ such that $(u_{n},v_{n})\rightarrow(u_{\varepsilon},v_{\varepsilon})$, $c_{\varepsilon}=I_{\varepsilon}(u_{\varepsilon},v_{\varepsilon})$ and $I_{\varepsilon}'(u_{\varepsilon},v_{\varepsilon})=0$. Moreover, since $\beta>\beta_{0}^\varepsilon$, from Lemma \ref{Lemma 3.6.}, we know that $u_{\varepsilon}\not\equiv0$ and $u_{\varepsilon}\not\equiv0$. Thus we prove $(u_{\varepsilon},v_{\varepsilon})\in\mathcal{N}_{\varepsilon}$. Finally, we prove that $u_{\varepsilon},v_{\varepsilon}>0$. Since $(|u_{\varepsilon}|,|v_{\varepsilon}|)\in\mathcal{N}_{\varepsilon}$ and $c_{\varepsilon}=I_{\varepsilon}(|u_{\varepsilon}|,|v_{\varepsilon}|)$, we deduce that $(|u_{\varepsilon}|,|v_{\varepsilon}|)$ is a nonnegative solution of (\ref{question3}). Using elliptic estimates and a Harnack inequality, we infer that $|u_{\varepsilon}|,|v_{\varepsilon}|>0$.  The proof is completed.
\qed

\s{Proof of Theorem \ref{Theorem 1.1.}}
\renewcommand{\theequation}{4.\arabic{equation}}
\begin{lemma}
\label{Lemma 4.1.}\cite{Schaftingen} {\it Let a function $V\in C(\mathbb{R}^{N},\mathbb{R})$ and $\gamma>0$. If $V\geq0$ on $\mathbb{R}^{N}$ with $\liminf_{|x|\rightarrow\infty}V(x)>0$,
and if for each $x\in\mathbb{R}^{N}$ such that
\begin{align*}
\liminf_{z\rightarrow x}\frac{V(z)}{|z-x|^{\gamma}}>0,
\end{align*}
then there exist $l\in\mathbb{N}$, $a_{1},\cdots,a_{l}\in\mathbb{R}^{N}$, $\mu>0$ and $\nu>0$ such that for each  $x\in\mathbb{R}^{N}$,
\begin{align*}
V(x)\geq\min\left\{\mu,\nu|x-a_{1}|^{\gamma},\cdots, \nu|x-a_{l}|^{\gamma}\right\}.
\end{align*}
}
\end{lemma}
\


\noindent Let us set
\begin{align*}
k=-\frac{\gamma}{2+\gamma}.
\end{align*}

\
\begin{proposition}
\label{Proposition 4.3.} {\it Under the assumptions of Theorem \ref{Theorem 1.1.}, there holds
\begin{align*}
\mathop{\limsup}\limits_{\varepsilon\rightarrow0}\frac{c_{\varepsilon}}{\varepsilon^{k(N-4)}}\leq m=\mathop{\inf}\limits_{i\in\mathcal{I}}\mathcal{E}_{W_{x_i}}.
\end{align*}
}
\end{proposition}
\

\noindent{\it Proof.} For $i\in \mathcal{I}$, $x_{i}\in\Omega$, let $W_{i}\in C(\mathbb{R}^{N})$ be the positive $\gamma$-homogeneous function such that
\begin{align}\label{15}
\mathop{\lim}\limits_{x\rightarrow x_{i}}\frac{a(x)-W_{i}(x-x_{i})}{|x-x_{i}|^\gamma}=0\ \ \mbox{and}\ \ \mathop{\lim}\limits_{x\rightarrow x_{i}}\frac{b(x)-W_{i}(x-x_{i})}{|x-x_{i}|^\gamma}=0,
\end{align}
and let $(\eta,\varphi)\in \mathbb{C}_{0}^{\infty}(\mathbb{R}^{N})\setminus\left\{(0,0)\right\}$, we define $\eta_{\varepsilon}(\mbox{resp.}\varphi_{\varepsilon}):\mathbb{R}^{N}\rightarrow\mathbb{R}$ by
\begin{align*}
(\eta_{\varepsilon}(x),\varphi_{\varepsilon}(x))=\varepsilon^{-k}\big(\eta(\frac{x-x_{i}/\varepsilon}{\varepsilon^{k}}),\varphi(\frac{x-x_{i}/\varepsilon}{\varepsilon^{k}})\big).
\end{align*}
Computing directly, we have
\begin{align*}
&\int_{\mathbb{R}^{N}}|\nabla \eta_{\varepsilon}|^{2}+\int_{\mathbb{R}^{N}}|\nabla \varphi_{\varepsilon}|^{2}=\varepsilon^{k(N-4)}\int_{\mathbb{R}^{N}}|\nabla \eta|^{2}+|\nabla \varphi|^{2};\\
&\int_{\mathbb{R}^{N}}a(\varepsilon x)|\eta_{\varepsilon}|^{2}+b(\varepsilon x)|\varphi_{\varepsilon}|^{2}=\varepsilon^{k(N-4)}\int_{\mathbb{R}^{N}}\frac{a(\varepsilon^{k+1}x+x_{i})}{\varepsilon^{(k+1)\gamma}}|\eta|^{2}
+\frac{b(\varepsilon^{k+1}x+x_{i})}{\varepsilon^{(k+1)\gamma}}|\varphi|^{2};\\
&\int_{\mathbb{R}^{N}}\mu_{1}\eta_{\varepsilon}^{4}+2\beta\eta_{\varepsilon}^{2}\varphi_{\varepsilon}^{2}+\mu_{2}\varphi_{\varepsilon}^{4}
=\varepsilon^{k(N-4)}\int_{\mathbb{R}^{N}}\mu_{1}\eta^{4}+2\beta\eta^{2}\varphi^{2}+\mu_{2}\varphi^{4}.
\end{align*}
Since the function $W_{i}$ is $\gamma$-homogeneous and satisfies (\ref{15}), we have
\begin{align}
\mathop{\lim}\limits_{\varepsilon\rightarrow0}\frac{a(\varepsilon^{k+1}x+x_{i})}{\varepsilon^{(k+1)\gamma}}
=\mathop{\lim}\limits_{\varepsilon\rightarrow0}\frac{W_{i}(\varepsilon^{k+1}x)}{\varepsilon^{(k+1)\gamma}}=W_{i}(x)\ \ \mbox{and}\ \ \mathop{\lim}\limits_{\varepsilon\rightarrow0}\frac{b(\varepsilon^{k+1}x+x_{i})}{\varepsilon^{(k+1)\gamma}}=W_{x_i}(x),
\end{align}
uniformly for $x\in supp(\eta)\cap supp(\varphi)$. Thus by Lebesgue's dominated convergence theorem, we have
\begin{align*}
\mathop{\lim}\limits_{\varepsilon\rightarrow0}\frac{1}{\varepsilon^{k(N-4)}}\int_{\mathbb{R}^{N}}a(\varepsilon x)|\eta_{\varepsilon}|^{2}+b(\varepsilon x)|\varphi_{\varepsilon}|^{2}
=\int_{\mathbb{R}^{N}}W_{i}|\eta|^{2}+W_{i}|\varphi|^{2}.
\end{align*}
Choosing $t_{\varepsilon}>0$ such that $\langle I_{\varepsilon}'(t_{\varepsilon}\eta_{\varepsilon},t_{\varepsilon}\varphi_{\varepsilon}),(t_{\varepsilon}\eta_{\varepsilon},t_{\varepsilon}\varphi_{\varepsilon})\rangle=0$, we observe that $t_{\varepsilon}\rightarrow t_{*}\in(0,+\infty)$ as $\varepsilon\rightarrow0$. Moreover, there holds
\begin{align*}
\int_{\mathbb{R}^{N}}|\nabla \eta|^{2}+W_{i}|\eta|^{2}+|\nabla \varphi|^{2}+W_{i}|\varphi|^{2}=t_{*}^2\int_{\mathbb{R}^{N}}\mu_{1}\eta^{4}+2\beta\eta^{2}\varphi^{2}+\mu_{2}\varphi^{4},
\end{align*}
which implies
\begin{align*}
\mathop{\limsup}\limits_{\varepsilon\rightarrow0}\frac{c_{\varepsilon}}{\varepsilon^{k(N-4)}}
\leq\mathop{\lim}\limits_{\varepsilon\rightarrow0}\frac{I_{\varepsilon}(t_{\varepsilon}\eta_{\varepsilon},t_{\varepsilon}\varphi_{\varepsilon})}{\varepsilon^{k(N-4)}}
=\mathcal{J}_{W_{i}}(t_{*}\eta,t_{*}\varphi)=\mathop{\sup}\limits_{t\in(0,+\infty)}\mathcal{J}_{W_{i}}(t\eta,t\varphi).
\end{align*}
Then by density of $\mathbb{C}_{0}^{\infty}(\mathbb{R}^{N})$ in the space $E_{W}:=H_W\times H_W$, by taking the infimum with respect to $\eta$ and $\varphi$, we have
\begin{align*}
\mathop{\limsup}\limits_{\varepsilon\rightarrow0}\frac{c_{\varepsilon}}{\varepsilon^{k(N-4)}}\leq\mathcal{E}_{W_{i}}.
\end{align*}
By the arbitrary choice of $i\in\mathcal{I}$, we obtain the conclusion.  The proof is completed.
\qed

\

\

\noindent{\bf Proof of  Theorem \ref{Theorem 1.1.}}
\bp Let $(u_{\varepsilon},v_{\varepsilon})$ be any positive ground state of system (\ref{question3}) obtained in Proposition \ref{Proposition 3.7.} and satisfying $I_{\varepsilon}(u_{\varepsilon},v_{\varepsilon})=c_{\varepsilon}$. Define
\begin{align*}
\omega_{\varepsilon}=\varepsilon^{k}u_{\varepsilon}(\varepsilon^{k}\cdot),\\
\phi_{\varepsilon}=\varepsilon^{k}v_{\varepsilon}(\varepsilon^{k}\cdot),
\end{align*}
one has
\begin{align}\label{18}
\begin{split}
&\int_{\mathbb{R}^{N}}|\nabla u_{\varepsilon}|^{2}=\varepsilon^{k(N-4)}\int_{\mathbb{R}^{N}}|\nabla \omega_{\varepsilon}|^{2};\\
&\int_{\mathbb{R}^{N}}a(\varepsilon x)u_{\varepsilon}^{2}=\varepsilon^{k(N-4)}\int_{\mathbb{R}^{N}}\e^{2k}a_{\varepsilon}(x)|\omega_{\varepsilon}|^{2};\\
&\int_{\mathbb{R}^{N}}\mu_{1}u_{\varepsilon}^{4}=\varepsilon^{k(N-4)}\int_{\mathbb{R}^{N}}\mu_{1}\omega_{\varepsilon}^{4};\\
&\int_{\mathbb{R}^{N}}\beta u_{\varepsilon}^{2}v_{\varepsilon}^{2}=\varepsilon^{k(N-4)}\int_{\mathbb{R}^{N}}\beta\omega_{\varepsilon}^{2}\phi_{\varepsilon}^{2}.\\
\end{split}
\end{align}
Similarly, one can get that
\begin{align}
-\triangle\omega_{\varepsilon}+\varepsilon^{2k}a_{\varepsilon}(x)\omega_{\varepsilon}=\mu_{1}\omega_{\varepsilon}^{3}+\beta\omega_{\varepsilon}\phi_{\varepsilon}^{2};\label{22}\\
-\triangle\phi_{\varepsilon}+\varepsilon^{2k}b_{\varepsilon}(x)\phi_{\varepsilon}=\mu_{2}\phi_{\varepsilon}^{3}+\beta\omega_{\varepsilon}^{2}\phi_{\varepsilon},\label{23}
\end{align}
where $a_{\varepsilon}(x)=a(\varepsilon^{k+1}x)$, $b_{\varepsilon}(x)=b(\varepsilon^{k+1}x)$. Consequently, there holds
\begin{align}\label{19}
\int_{\mathbb{R}^{N}}|\nabla \omega_{\varepsilon}|^{2}+\varepsilon^{2k}a_{\varepsilon}(x)\omega_{\varepsilon}^{2}+|\nabla \phi_{\varepsilon}|^{2}+\varepsilon^{2k}b_{\varepsilon}(x)\phi_{\varepsilon}^{2}=\int_{\mathbb{R}^{N}}\mu_{1}\omega_{\varepsilon}^{4}+2\beta\omega_{\varepsilon}^{2}\phi_{\varepsilon}^{2}+\mu_{2}\phi_{\varepsilon}^{4},
\end{align}
By Proposition \ref{Proposition 4.3.}, we deduce that, for small $\varepsilon>0$,
\begin{align}\label{20}
\int_{\mathbb{R}^{N}}|\nabla \omega_{\varepsilon}|^{2}+\varepsilon^{2k}a_{\varepsilon}(x)\omega_{\varepsilon}^{2}+|\nabla \phi_{\varepsilon}|^{2}+\varepsilon^{2k}b_{\varepsilon}(x)\phi_{\varepsilon}^{2}\leq4m.
\end{align}

{\bf Step 1.} Now we claim that there exists a sequence $\{x_{\varepsilon}\}\in\mathbb{R}^{N}$ such that
\begin{align*}
\liminf_{\varepsilon\rightarrow0}\int_{B_{1}(x_{\varepsilon})}\omega_{\varepsilon}^{4}+\phi_{\varepsilon}^{4}>0.
\end{align*}
It suffices to show that
\begin{align}\label{30}
\lim_{\e\rightarrow0}\sup_{x\in\R^{N}}\varepsilon^{k(4-N)}\int_{B_{\varepsilon^{k}(x)}}u_{\e}^{4}+v_{\e}^{4}>0.
\end{align}
We adopt some ideas in \cite{Schaftingen} to show that (\ref{30}) holds. By the scaled version of the classical Sobolev embedding theorem, for every $x\in\R^{N}$
\begin{align*}
\Bigg(\int_{B_{\varepsilon^{k}(x)}}u_{\varepsilon}^{4}\Bigg)^{\frac{1}{2}}\leq C_{1}\varepsilon^{-\frac{kN}{2}}\int_{B_{\varepsilon^{k}(x)}}\varepsilon^{2k}|\nabla u_{\varepsilon}|^{2}+u_{\varepsilon}^{2},
\end{align*}
where $C_{1}$ is independent of $x\in\R^{N}$ and $\varepsilon$. By \cite[Lemma 3.4]{Schaftingen}, for some constant $C_{2}$ independent of $x\in\R^{N}$ and $\varepsilon$ small enough, we obtain
\begin{align*}
\Bigg(\int_{B_{\varepsilon^{k}(x)}}u_{\varepsilon}^{4}\Bigg)^{\frac{1}{2}}\leq C_{2}\varepsilon^{-\frac{kN}{2}+2k}\int_{B_{\varepsilon^{k}(x)}}|\nabla u_{\e}|^{2}+a(\varepsilon x)u_{\e}^{2},
\end{align*}
and then
\begin{align}\label{31}
\Bigg(\int_{B_{\varepsilon^{k}(x)}}u_{\e}^{4}\Bigg)\leq C_{2}\varepsilon^{-\frac{kN}{2}+2k}\Bigg(\int_{B_{\varepsilon^{k}(x)}}u_{\varepsilon}^{4}\Bigg)^{\frac{1}{2}}\int_{B_{\varepsilon^{k}(x)}}|\nabla u_{\e}|^{2}+a(\varepsilon x)u_{\e}^{2}.
\end{align}
By integration both sides on (\ref{31}) with respect to $x$ over $\R^{N}$,
\begin{align*}
\int_{\R^{N}}u_{\e}^{4}\leq C_{2}\varepsilon^{-\frac{kN}{2}+2k}\Bigg(\sup_{x\in\R^{N}}\int_{B_{\varepsilon^{k}(x)}}u_{\varepsilon}^{4}\Bigg)^{\frac{1}{2}}\int_{\R^{N}}|\nabla u_{\varepsilon}|^{2}+a(\varepsilon x)u_{\varepsilon}^{2}.
\end{align*}
Similarly ,
\begin{align*}
\int_{\R^{N}}v_{\e}^{4}\leq C_{2}\varepsilon^{-\frac{kN}{2}+2k}\Bigg(\sup_{x\in\R^{N}}\int_{B_{\varepsilon^{k}(x)}}v_{\varepsilon}^{4}\Bigg)^{\frac{1}{2}}\int_{\R^{N}}|\nabla v_{\e}|^{2}+b(\varepsilon x)v_{\e}^{2}.
\end{align*}
Denote $C^{*}=\max\Big\{\Big(\sup_{x\in\R^{N}}\int_{B_{\varepsilon^{k}(x)}}u_{\e}^{4}\Big)^{\frac{1}{2}},\Big(\sup_{x\in\R^{N}}\int_{B_{\varepsilon^{k}(x)}}v_{\e}^{4}\Big)^{\frac{1}{2}}\Big\}$,
\begin{align*}
\int_{\R^{N}}u_{\e}^{4}+v_{\e}^{4}\leq C^{*}C_{2}\varepsilon^{-\frac{kN}{2}+2k}\int_{\R^{N}}|\nabla u_{\e}|^{2}+a(\varepsilon x)u_{\e}^{2}+|\nabla v_{\e}|^{2}+b(\varepsilon x)v_{\e}^{2}
\end{align*}
can be obtained.
Since $(u_{\varepsilon},v_{\varepsilon})$ is a solution of problem (\ref{question3}) satisfing
\begin{align*}
&0<\int_{\R^{N}}|\nabla u_{\e}|^{2}+a(\varepsilon x)u_{\e}^{2}+|\nabla v_{\e}|^{2}+b(\varepsilon x)v_{\e}^{2}\\
&=\int_{\R^{N}}\mu_{1}u_{\e}^{4}+2\beta u_{\e}^{2}v_{\e}^{2}+\mu_{2}v_{\e}^{4}\leq (\mu_1+\mu_2+\beta)\int_{\R^{N}}u_{\e}^{4}+v_{\e}^{4}\\
&\leq C^{*}C_2(\mu_1+\mu_2+\beta)\varepsilon^{-\frac{kN}{2}+2k}\int_{\R^{N}}|\nabla u_{\e}|^{2}+a(\varepsilon x)u_{\e}^{2}+|\nabla v_{\e}|^{2}+b(\varepsilon x)v_{\e}^{2},
\end{align*}
one get that
$$
(C^{*})^2\varepsilon^{-\frac{kN}{2}+2k}\ge\left(C_2(\mu_1+\mu_2+\beta)\right)^{-2}.
$$
It follows that
\begin{align*}
\liminf_{\varepsilon\rightarrow0}\sup_{x\in\R^{N}}\varepsilon^{k(4-N)}\int_{B_{\varepsilon^{k}(x)}}u_{\e}^{4}>0,
\end{align*}
or
\begin{align*}
\liminf_{\varepsilon\rightarrow0}\sup_{x\in\R^{N}}\varepsilon^{k(4-N)}\int_{B_{\varepsilon^{k}(x)}}v_{\e}^{4}>0.
\end{align*}
As a consequence, \eqref{30} holds.

{\bf Step 2.}
We claim that
\begin{align}\label{CC}
\liminf_{\e\rightarrow0}\int_{B_{1}(x_{\varepsilon})}\omega_{\varepsilon}^{4}>0\ \ \mbox{and}\ \  \liminf_{\e\rightarrow0}\int_{B_{1}(x_{\varepsilon})}\phi_{\varepsilon}^{4}>0.
\end{align}
We argue by contradiction. Suppose that
\begin{align*}
\limsup_{\e\rightarrow0}\sup_{x\in\R^{N}}\int_{B_{1}(x)}\omega_{\e}^{4}=0,
\end{align*}
by the Lions Concentration-Compactness Principle\cite{PL}, we deduce that $\omega_{\e}\rightarrow0$ in $L^{p}(\R^{N})$ for $2<p<2^{*}$. Multiplying \eqref{22} and \eqref{23} by $\omega_{\e}$ and $\phi_{\e}$ respectively and integrating both sides with respect to $x$ over $\R^{N}$, we obtain
\begin{align*}
\int_{\mathbb{R}^{N}}|\nabla \omega_{\varepsilon}|^{2}+\varepsilon^{2k}a_{\varepsilon}(x)\omega_{\varepsilon}^{2}=o_{\e}(1),
\end{align*}
\begin{align}\label{32}
\int_{\mathbb{R}^{N}}|\nabla \phi_{\varepsilon}|^{2}+\varepsilon^{2k}b_{\varepsilon}(x)\phi_{\varepsilon}^{2}=\int_{\mathbb{R}^{N}}\mu_{2}\phi_{\varepsilon}^{4}+o_{\e}(1).
\end{align}
Therefore, there exists $t_{\e}>0$ such that
\begin{align}\label{33}
\int_{\mathbb{R}^{N}}|\nabla \phi_{\varepsilon}|^{2}+\varepsilon^{2k}b_{\varepsilon}(x)\phi_{\varepsilon}^{2}=t_{\e}^{2}\int_{\mathbb{R}^{N}}\mu_{2}\phi_{\varepsilon}^{4}.
\end{align}
Combing the fact that $\liminf_{\e\rightarrow0}\int_{B_{1}(x_{\varepsilon})}\phi_{\varepsilon}^{4}>0$, it follows that $t_{\e}\rightarrow1$ by \eqref{32} and \eqref{33}, as $\e\rightarrow0$. Hence
\begin{align*}
\lim_{\varepsilon\rightarrow0}\varepsilon^{k(4-N)}I_{\e}(u_{\e},v_{\e})&=\lim_{\varepsilon\rightarrow0}\left[\frac{1}{2}\int_{\mathbb{R}^{N}}(|\nabla \phi_{\varepsilon}|^{2}+\varepsilon^{2k}b_{\varepsilon}(x)\phi_{\varepsilon}^{2})-\frac{\mu_{2}}{4}\int_{\mathbb{R}^{N}}\phi_{\varepsilon}^{4}\right]\\
&=\lim_{\varepsilon\rightarrow0}\left[\frac{t_\varepsilon^2}{2}\int_{\mathbb{R}^{N}}(|\nabla \phi_{\varepsilon}|^{2}+\varepsilon^{2k}b_{\varepsilon}(x)\phi_{\varepsilon}^{2})-\frac{\mu_{2}t_\varepsilon^4}{4}\int_{\mathbb{R}^{N}}\phi_{\varepsilon}^{4}\right]\\
&=\lim_{\varepsilon\rightarrow0}\varepsilon^{k(4-N)}I_{\e}(0,t_\varepsilon v_{\e}).
\end{align*}
Noting that $I_{\e}(0,t_\varepsilon v_{\e})\ge I_{\e}(0,V_{1\e}(\varepsilon\cdot))$, one has
\begin{equation}\label{energy}
\lim_{\varepsilon\rightarrow0}\varepsilon^{k(4-N)}c_{\e}\ge\lim_{\varepsilon\rightarrow0}\varepsilon^{k(4-N)}I_{\e}(0,V_{1\e}(\varepsilon\cdot)).
\end{equation}
By rescaling $U_{1\e}$ and $V_{1\e}$, $f_{\e}$ obtains its minimum at $(s_\varepsilon,t_\varepsilon)$, where
$$
s_\varepsilon=\mu_1 \varepsilon^{-\frac{8\gamma}{\gamma+2}}\int_{\R^{N}}U_{1\e}^{4}(\e^{\frac{2}{\gamma+2}}x)\left[\beta\int_{\R^{N}}U_{1\e}^{2}(\e^{\frac{2}{\gamma+2}}x)V_{1\e}^2(\e^{\frac{2}{\gamma+2}}x)
-\mu_2\int_{\R^{N}}V_{1\e}^{4}(\e^{\frac{2}{\gamma+2}}x)\right]
$$
and
$$
t_\varepsilon=\mu_2 \varepsilon^{-\frac{8\gamma}{\gamma+2}}\int_{\R^{N}}V_{1\e}^{4}(\e^{\frac{2}{\gamma+2}}x)\left[\beta\int_{\R^{N}}U_{1\e}^{2}(\e^{\frac{2}{\gamma+2}}x)V_{1\e}^2(\e^{\frac{2}{\gamma+2}}x)
-\mu_1\int_{\R^{N}}U_{1\e}^{4}(\e^{\frac{2}{\gamma+2}}x)\right].
$$
It follows from Remark \ref{remark} that $s_\varepsilon\rightarrow s_0$, $t_\varepsilon\rightarrow t_0$ as $\varepsilon\rightarrow0$,
where
$$
s_0=\mu_1 a_0[\beta c_0-\mu_2b_0],\,\ t_0=\mu_2 b_0[\beta c_0-\mu_1a_0],
$$
and
$$
a_0=\int_{\R^{N}}\omega_{\gamma,\mu_1}^{4},\,\,b_0=\int_{\R^{N}}\omega_{\gamma,\mu_2}^{4},\,\,c_0=\int_{\R^{N}}\omega_{\gamma,\mu_1}^{2}\omega_{\gamma,\mu_2}^{2}.
$$
Meanwhile,
\begin{align*}
\varepsilon^{k(4-N)}f_{\e}(s,t)=f_{0}(s,t)+o_{\e}(1)\ \mbox{uniformly for $(s,t)\in \Gamma$, as}\ \e\rightarrow0,
\end{align*}
where
\begin{align*}
f_{0}(s,t)=\frac{(s\mu_1a_0+t\mu_2 b_0)^{2}}{4(\mu_{1}s^{2}a_0+2\beta stc_0+\mu_{2}t^{2}b_0)}.
\end{align*}
It follows that
\begin{align*}
f_{0}(s,0)=\frac{\mu_1a_0}{4}=\lim_{\varepsilon\rightarrow0}\varepsilon^{k(4-N)}I_{\e}(U_{1\e}(\varepsilon\cdot),0),\\
f_{0}(0,t)=\frac{\mu_2b_0}{4}=\lim_{\varepsilon\rightarrow0}\varepsilon^{k(4-N)}I_{\e}(0,V_{1\e}(\varepsilon\cdot)).
\end{align*}
Moreover, for $\beta>\beta_0$,
$$
f_0(s_0,t_0)=\min_{(s,t)\in\Gamma}f_{0}(s,t)<\min\{f_0(s,0),f_0(0,t)\}.
$$
Then we observe that
$$
\lim_{\varepsilon\rightarrow0}\varepsilon^{k(4-N)}f_{\e}(s_\varepsilon,t_\varepsilon)<\lim_{\varepsilon\rightarrow0}
\varepsilon^{k(4-N)}\min\left\{I_{\e}(U_{1\e}(\varepsilon\cdot),0),I_{\e}(V_{1\e}(\varepsilon\cdot),0)\right\}.
$$
By Lemma \ref{Lemma 3.6.},
$$
\limsup_{\varepsilon\rightarrow0}\varepsilon^{k(4-N)}c_\varepsilon<\lim_{\varepsilon\rightarrow0}
\varepsilon^{k(4-N)}\min\left\{I_{\e}(U_{1\e}(\varepsilon\cdot),0),I_{\e}(V_{1\e}(\varepsilon\cdot),0)\right\},
$$
which contradicts \eqref{energy}.

{\bf Step 3.}
We claim that $\varepsilon^{k+1}x_{\varepsilon}\rightarrow x_{i}$ for some $i\in\mathcal{I}$. By {\bf Step 1-2}, for some $\tilde{\omega},\tilde{\phi}\in H_{loc}^1(\mathbb{R}^N)\setminus\{0\}$, such that, up to a subsequence, $\omega_\varepsilon(\cdot+x_\e)\rightharpoonup \tilde{\omega}$ and $\phi_\varepsilon(\cdot+x_\e)\rightharpoonup \tilde{\phi}$ weakly in $H_{loc}^1(\mathbb{R}^N)$ as $\varepsilon\rightarrow0$. Then if $\varepsilon^{k+1}x_{\varepsilon}\rightarrow \infty$ as $\varepsilon\rightarrow0$, by Fatou's Lemma,
$$
\liminf_{\varepsilon\rightarrow\infty}\int_{\mathbb{R}^{N}}a_{\varepsilon}(x)\omega_{\varepsilon}^{2}+b_{\varepsilon}(x)\phi_{\varepsilon}^{2}
\ge\int_{\mathbb{R}^{N}}a_\infty\tilde{\omega}^2+b_\infty\tilde{\phi}^{2}.
$$
It follows from \eqref{20} that
$$
\int_{\mathbb{R}^{N}}a_\infty\tilde{\omega}^2+b_\infty\tilde{\phi}^{2}=0,
$$
which implies that $a_\infty=b_\infty=0$. This is a contradiction. As a consequence, up to a subsequence, $\varepsilon^{k+1}x_{\varepsilon}\rightarrow x_\ast$ for some $x_\ast\in\mathbb{R}^N$ as $\varepsilon\rightarrow0$. Then Fatou's Lemma tells us that
$$
\int_{\mathbb{R}^{N}}a(x_\ast)\tilde{\omega}^2+b(x_\ast)\tilde{\phi}^{2}=0
$$
and then $x_\ast\in\mathcal{A}\bigcap\mathcal{B}$

Now, we claim that $x_\ast=x_{i}$ for some $i\in\mathcal{I}$. By the assumption of Theorem \ref{Theorem 1.1.}, there exists some positive $\gamma$-homogeneous function $W_{\ast}\in C(\mathbb{R}^{N})$, such that
$$
\mathop{\lim}\limits_{x\rightarrow x_\ast}\frac{a(x)-W_{\ast}(x-x_\ast)}{|x-x_\ast|^{\gamma}}=0\,\mbox{or}\,\,\infty.
$$
For any $x\not=x_\ast$, let
$$
\lambda(x-x_\ast)=\frac{a(x)-W_{\ast}(x-x_\ast)}{|x-x_\ast|^{\gamma}},\,\,\mu(x-x_\ast)=\frac{b(x)-W_{\ast}(x-x_\ast)}{|x-x_\ast|^{\gamma}},
$$
then to conclude the claim, we just need to show that $\lim_{x\rightarrow 0}\lambda(x)=0$, which implies that $\lim_{x\rightarrow 0}\mu(x)=0$.

According to (\ref{20}),
\begin{align*}
\limsup_{\varepsilon\rightarrow0}\mathop{\inf}\limits_{x\in B_{1}(x_{\varepsilon})}\left[\frac{a_{\varepsilon}(x)}{\varepsilon^{(k+1)\gamma}}+\frac{b_{\varepsilon}(x)}{\varepsilon^{(k+1)\gamma}}\right]<\infty.
\end{align*}
Thus there exists $y_\varepsilon\in B_{1}(x_{\varepsilon})$ such that
$$
\limsup_{\varepsilon\rightarrow0}\left[\frac{a(\varepsilon^{k+1}y_\varepsilon)}{\varepsilon^{(k+1)\gamma}}+\frac{b(\varepsilon^{k+1}y_\varepsilon)}{\varepsilon^{(k+1)\gamma}}\right]<\infty
$$
and $\varepsilon^{k+1}y_\varepsilon\rightarrow x_\ast$ as $\varepsilon\rightarrow0$.
For $\varepsilon$ small, by Lemma \ref{Lemma 4.1.}, there holds $a(\varepsilon^{k+1}y_\varepsilon)\ge c|\varepsilon^{k+1}y_\varepsilon-x_\ast|^\gamma$ and $b(\varepsilon^{k+1}y_\varepsilon)\ge c|\varepsilon^{k+1}y_\varepsilon-x_\ast|^\gamma$ for some $c>0$. Then for $\varepsilon$ small and some $C>0$,
\begin{align*}
|\varepsilon^{k+1}y_\varepsilon-x_\ast|\le C\varepsilon^{k+1},
\end{align*}
from which it follows that
\begin{align*}
|x_{\varepsilon}-\frac{x_{*}}{\varepsilon^{k+1}}|\leq C+1.
\end{align*}
So for $R>0$ large enough, we have $B_1(x_\varepsilon)\subset B_{R}(\frac{x_{*}}{\varepsilon^{k+1}})$. Due to \eqref{CC}, for some $\omega,\phi\in H_{loc}^1(\mathbb{R}^N)\setminus\{0\}$, such that, up to a subsequence, $\tilde{\omega}_\varepsilon:=\omega_\varepsilon(\cdot+\frac{x_{*}}{\varepsilon^{k+1}})\rightharpoonup \omega$ and $\tilde{\phi}_\varepsilon:=\phi_\varepsilon(\cdot+\frac{x_{*}}{\varepsilon^{k+1}})\rightharpoonup \phi$ weakly in $H_{loc}^1(\mathbb{R}^N)$ as $\varepsilon\rightarrow0$.

If $\lim_{x\rightarrow 0}\lambda(x)=+\infty$, then for any $x\in\mathbb{R}^N$,
$$
\varepsilon^{2k}a(\varepsilon^{k+1}x+x_\ast)=W_{\ast}(x)+\lambda(\varepsilon^{k+1}x)|x|^\gamma
$$
and
$$
\varepsilon^{2k}b(\varepsilon^{k+1}x+x_\ast)=W_\ast(x)+\mu(\varepsilon^{k+1}x)|x|^\gamma.
$$
Thanks to \eqref{20} and $\mu(\cdot)$ is bounded from below in $\mathbb{R}^N$, by Rellich's theorem,
\begin{align*}
4m&\ge\liminf_{\varepsilon\rightarrow0}\int_{B_R(\frac{x_{*}}{\varepsilon^{k+1}})}\varepsilon^{2k}a_{\varepsilon}(x)\omega_{\varepsilon}^{2}+\varepsilon^{2k}b_{\varepsilon}(x)\phi_{\varepsilon}^{2}\\
&=\liminf_{\varepsilon\rightarrow0}\int_{B_R(0)}\varepsilon^{2k}a(\varepsilon^{k+1}x+x_\ast)\tilde{\omega}_{\varepsilon}^{2}(x)
+\varepsilon^{2k}b(\varepsilon^{k+1}x+x_\ast)\tilde{\phi}_{\varepsilon}^{2}(x)\\
&=\liminf_{\varepsilon\rightarrow0}\int_{B_R(0)}W_\ast(x)(\tilde{\omega}_{\varepsilon}^{2}(x)+\tilde{\phi}_{\varepsilon}^{2}(x))
+|x|^\gamma(\lambda(\varepsilon^{k+1}x)\tilde{\omega}_{\varepsilon}^{2}(x)+\mu(\varepsilon^{k+1}x)\tilde{\phi}_{\varepsilon}^{2}(x))\\
&\ge\liminf_{\varepsilon\rightarrow0}\int_{B_R(0)}|x|^\gamma(\lambda(\varepsilon^{k+1}x)
\tilde{\omega}_{\varepsilon}^{2}(x)+\mu(\varepsilon^{k+1}x)\tilde{\phi}_{\varepsilon}^{2}(x))=\infty,
\end{align*}
which is a contradiction. Thus,
$$
\lim_{x\rightarrow x_\ast}\frac{a(x)-W_\ast(x-x_\ast)}{|x-x_\ast|^{\gamma}}=\lim_{x\rightarrow x_\ast}\frac{b(x)-W_\ast(x-x_\ast)}{|x-x_\ast|^{\gamma}}=0.
$$
By the assumption of Theorem \ref{Theorem 1.1.}, $x_\ast=x_i$ for some $i\in\mathcal{I}$.

{\bf Step 4.} It is easy to know
\begin{align*}
\int_{\mathbb{R}^{N}}|\nabla\tilde{\omega}_{\varepsilon}|^{2}+\e^{2k}a(\varepsilon^{k+1}x+x_\ast)\tilde{\omega}_{\varepsilon}^{2}
+|\nabla\tilde{\phi}_{\varepsilon}|^{2}+\e^{2k}b(\varepsilon^{k+1}x+x_\ast)\tilde{\phi}_{\varepsilon}^{2}
=\int_{\mathbb{R}^{N}}\mu_{1}\tilde{\omega}_{\varepsilon}^{4}+2\beta\tilde{\omega}_{\varepsilon}^{2}\tilde{\phi}_{\varepsilon}^{2}+\mu_{2}\tilde{\phi}_{\varepsilon}^{4},
\end{align*}
Next, we prove $(\omega,\phi)$ is a critical point of $\mathcal{J}_{W_{\ast}}$. By Fatou's Lemma and the lower semi-continuity of the norm,
\begin{align*}
&\int_{B_R(0)}|\nabla \omega|^{2}+W_{\ast}\omega^{2}+|\nabla \phi|^{2}+W_{\ast}\phi^{2}\\
&\leq\liminf_{\varepsilon\rightarrow0}\int_{B_R(0)}|\nabla\tilde{\omega}_{\varepsilon}|^{2}+\e^{2k}a(\varepsilon^{k+1}x+x_\ast)\tilde{\omega}_{\varepsilon}^{2}
+|\nabla\tilde{\phi}_{\varepsilon}|^{2}+\e^{2k}b(\varepsilon^{k+1}x+x_\ast)\tilde{\phi}_{\varepsilon}^{2}\\
&\leq4m.
\end{align*}
Since $R$ is arbitrary, we get that
\begin{equation}\label{semi}
\int_{\mathbb{R}^N}|\nabla \omega|^{2}+W_{\ast}\omega^{2}+|\nabla \phi|^{2}+W_{\ast}\phi^{2}\le 4m
\end{equation}
and then
$\omega,\phi\in H^1(\mathbb{R}^N)$.
For any $\xi,\eta\in C_{0}^{\infty}(\mathbb{R}^{N})$, Multiplying equation (\ref{22}) and (\ref{23}) by $\xi\left(\cdot-\frac{x_{*}}{\varepsilon^{k+1}}\right)$ and $\eta\left(\cdot-\frac{x_{*}}{\varepsilon^{k+1}}\right)$ and making $\varepsilon\rightarrow0$,
\begin{align*}
\int_{\mathbb{R}^{N}}\nabla \omega\nabla\xi+W_{*}\omega\xi+\nabla \phi\nabla\eta+W_{*}\phi\eta=\int_{\mathbb{R}^{N}}\mu_{1}\omega^{3}\xi+\beta\phi^{2}\omega\xi+\beta \omega^{2}\phi\eta+\mu_{2}\phi^{3}\eta,
\end{align*}
which implies $\mathcal{J}_{W_{*}}'(\omega,\phi)=0$ in $(H_{W_\ast}\times H_{W_\ast})^\ast$ and $\mathcal{J}_{W_{*}}(\omega,\phi)\ge\mathcal{E}_{W_\ast}$. Due to \eqref{semi}, $\mathcal{J}_{W_{*}}(\omega,\phi)\le m=\mathop{\inf}\limits_{i\in\mathcal{I}}\mathcal{E}_{W_{i}}$. Therefore,
$$
\mathcal{J}_{W_{*}}(\omega,\phi)=\mathop{\inf}\limits_{i\in\mathcal{I}}\mathcal{E}_{W_{i}}.
$$
and $(\tilde{\omega}_{\varepsilon},\tilde \phi_{\varepsilon})\rightarrow(\omega,\phi)$ in $\mathbb{H}_{loc}$ as $\varepsilon\rightarrow0$.  The proof is completed.
\ep

\s{Proof of Theorem \ref{Theorem 1.2.}}
\renewcommand{\theequation}{5.\arabic{equation}}

This last section is devoted to the proof of Theorem \ref{Theorem 1.2.}, which covers the case that the potential $a$ and $b$ vanish on the closure of some smooth bounded open set.

\

\noindent{\bf Proof of Theorem \ref{Theorem 1.2.}.} \
\bp
{\bf Step 1.} Similarly to Theorem \ref{Theorem 1.1.}, for sufficiently small $\varepsilon>0$ and $\beta>\beta_0^\e$, we prove that system (\ref{question2}) has a positive ground state solution $(U_{\varepsilon},V_{\varepsilon})$. Define
$$
E_\varepsilon:=\left\{(u,v)\in \mathbb{H}_{loc}: \int_{\mathbb{R}^N}\varepsilon^2|\nabla u|^2+a(x)u^2<\infty,\,\,\int_{\mathbb{R}^N}\varepsilon^2|\nabla v|^2+b(x)v^2<\infty \right\}
$$
endowed with norm
$$
\|(u,v)\|_\varepsilon=\left(\int_{\mathbb{R}^N}\varepsilon^2|\nabla u|^2+a(x)u^2+\int_{\mathbb{R}^N}\varepsilon^2|\nabla v|^2+b(x)v^2\right)^{\frac12}.
$$
Define the auxiliary functional $J_{\varepsilon}\in C^{1}(E_\varepsilon)$ by for each $(w,z)\in E_\varepsilon$,
\begin{align*}
J_{\varepsilon}(w,z)=\frac{1}{2}\int_{\mathbb{R}^{N}}|\nabla w|^{2}+\frac{a(x)}{\varepsilon^{2}}w^{2}+|\nabla z|^{2}+\frac{b(x)}{\varepsilon^{2}}z^{2}-\frac{1}{4}F(w,z).
\end{align*}
Set
\begin{align*}
d_\varepsilon:=\inf\{J_{\varepsilon}(w,z):\ (w,z)\in E_\varepsilon\setminus\{(0,0)\},\langle J'_{\varepsilon}(w,z),(w,z)\rangle=0\}.
\end{align*}
For any $(u,v)\in E_\ast=H_0^1(int(\Omega))\times H_0^1(int(\Omega))$, set
\begin{align*}
J_{\ast}(w,z)=\frac{1}{2}\int_{\Omega}|\nabla w|^{2}+|\nabla z|^{2}-\frac{1}{4}\int_{\Omega}\mu_{1}w^{4}+2\beta w^{2}z^{2}+\mu_{2}z^{4}.
\end{align*}
Similarly to Lemma \ref{Lemma 3.6.}, for sufficiently small $\varepsilon>0$ and $\beta>\beta_0^\e$,
\begin{align}\label{erergy-es}
d_{\varepsilon}=\varepsilon^{N-4}c_\varepsilon<\varepsilon^{N-4}\min\{I_{\varepsilon}(U_{1\varepsilon}(\e\cdot),0),I_{\varepsilon}(0,V_{1\varepsilon}(\e\cdot)\}.
\end{align}
Moreover, for sufficiently small $\e$, problem (\ref{question2}) admits a ground state solution $(U_{\e},V_{\e})$, which is fully nontrivial.
By virtue of the Nehari manifold method, one can show that for any $\beta>0$, the following problem
\begin{align}\label{limitsystem}
\begin{cases}
-\Delta w=\mu_{1}w^{3}+\beta z^{2}w, \   x\in\mbox{int}(\Omega),\\
-\Delta z=\mu_{2}z^{3}+\beta w^{2}z, \  x\in\mbox{int}(\Omega),\\
w, z\in H_{0}^{1}(\mbox{int}(\Omega))
\end{cases}
\end{align}
admits a ground state solution. Moreover, set
\begin{align*}
0<c_{\ast}=\inf\{J_{\ast}(w,z):\ (w,z)\in E_\ast\setminus\{(0,0)\},\langle J_{\ast}'(w,z),(w,z)\rangle=0\},
\end{align*}
it is easy to check that $d_\e=J_{\varepsilon}(\varepsilon^{-1}U_{\varepsilon},\varepsilon^{-1}V_{\varepsilon})\le c_\ast\le\tilde{c}_{\ast}$, where $\tilde{c}_{\ast}$ is the least energy of the following problem
\begin{align*}
\begin{cases}
-\Delta w=\beta z^{2}w, \   x\in\mbox{int}(\Omega),\\
-\Delta z=\beta w^{2}z, \  x\in\mbox{int}(\Omega),\\
w, z\in H_{0}^{1}(\mbox{int}(\Omega)).
\end{cases}
\end{align*}
\vskip0.1in
{\bf Step 2.} Let $w_{\varepsilon}=\varepsilon^{-1}U_{\varepsilon}$ and $z_{\varepsilon}=\varepsilon^{-1}V_{\varepsilon}$.
Noting that for any $\varepsilon>0$,
\begin{align}\label{34}
\int_{\mathbb{R}^{N}}|\nabla w_{\varepsilon}|^{2}+\frac{a(x)}{\varepsilon^{2}}w_{\varepsilon}^{2}+|\nabla z_{\varepsilon}|^{2}+\frac{b(x)}{\varepsilon^{2}}z_{\varepsilon}^{2}\leq4c_{\ast}.
\end{align}
It follows from \cite[Lemma 2.1]{Schaftingen} that there exists $C>0$(independent of $\varepsilon$) such that
\begin{align*}
\int_{\R^{N}}(|\na w_\e|^2+w_\e^2)\leq C\int_{\R^{N}}(|\na w_\e|^2+a(x)w_\e^2),
\end{align*}
which implies that $\int_{\R^{N}}(|\na w_\e|^2+w_\e^2)\leq 4Cc_{\ast}$ if $\varepsilon\le1$.
Similar result holds for $z_{\varepsilon}$. So, $(w_{\varepsilon},z_{\varepsilon})$ is bounded in $\mathbb{H}$ if $\varepsilon\le1$ and for some $(w,z)\in \mathbb{H}$, $(w_{\varepsilon},z_{\varepsilon})\rightharpoonup(w,z)$ in $\mathbb{H}$  and  $(w_{\varepsilon},z_{\varepsilon})\rightarrow(w,z)$ in $\mathbb{L}_{loc}^{2}(\mathbb{R}^{N})$ as $\varepsilon\rightarrow0$.

For $R>0$ large enough such that $a(x)\ge \frac{a_\infty}{2}$ and $b(x)\ge \frac{b_\infty}{2}$ for $|x|\ge R$. By \eqref{34},
\begin{align*}
\frac{1}{2}\min\{a_\infty,b_\infty\}\int_{\R^{N}\setminus B_R(0)}w_{\varepsilon}^{2}+z_{\e}^{2}\leq \int_{\R^{N}\setminus B_R(0)}a(x)w_{\varepsilon}^{2}+b(x)z_{\varepsilon}^{2}\le4c_\ast\varepsilon^2\rightarrow0,\ \ \ \varepsilon\rightarrow0.
\end{align*}
So Rellich's theorem tells us that $w_{\varepsilon}\rightarrow w$ and $z_{\varepsilon}\rightarrow z$ strongly in $L^{2}(\mathbb{R}^{N})$. On the other hand, due to
$J_{\varepsilon}'(w_{\varepsilon},z_{\varepsilon})=0$, for $\varepsilon>0$ small and some $C,C_1$(independent of $\varepsilon$),
\begin{align*}
\|(w_{\varepsilon},z_{\varepsilon})\|_{4}^{2}\leq C\|(w_{\varepsilon},z_{\varepsilon})\|_{\mathbb{H}}^{2}\leq C_{1}\|(w_{\varepsilon},z_{\varepsilon})\|_{4}^{4},
\end{align*}
which shows that $\|(w,z)\|_{4}>0$ and $(w,z)\neq(0,0)$. \

We claim that $w,z\not\equiv0$. Otherwise, without loss of generalization, if $z\equiv0$, then $z_{\varepsilon}\rightarrow 0$ strongly in $L^{2}(\mathbb{R}^{N})$ as $\e\rightarrow0$. Noting that $\{z_\e\}$ is bounded in $H^1(\mathbb{R}^N)$, by the Gagliardo-Nirenberg-Sobolev inequality, $z_\e\rightarrow0$ strongly in $L^{4}(\mathbb{R}^{N})$ as $\e\rightarrow0$ and
$$
\lim_{\varepsilon\rightarrow0}\int_{\mathbb{R}^{N}}|\nabla z_{\varepsilon}|^{2}+\frac{b(x)}{\varepsilon^{2}}z_{\varepsilon}^{2}=0.
$$
Then similarly to Section 4,
$$
\liminf_{\varepsilon\rightarrow0}d_\e\ge\liminf_{\varepsilon\rightarrow0}\e^{N-4}I_\e(U_{1\e}(\e \cdot),0),
$$
which gives a contradiction by using \eqref{erergy-es} similarly as above.

\vskip0.1in
{\bf Step 3.} We show that $w,z\in H_0^1(\mbox{int}(\Omega))$ for $\mu_1,\mu_2>0$ small. By \eqref{34}, we have
$$
\lim_{\varepsilon\rightarrow0}\int_{\mathbb{R}^N\setminus\Sigma}a(x)w^2+b(x)z^2=\lim_{\varepsilon\rightarrow0}\int_{\mathbb{R}^N\setminus\Sigma}a(x)w_\e^2+b(x)z_\e^2=0,
$$
which implies that $w=z=0$ in $\mathbb{R}^N\setminus\Sigma$. Similarly,
$$
\lim_{\varepsilon\rightarrow0}\int_{\Sigma\setminus\Omega}a(x)w^2+b(x)z^2=0
$$
which implies that $w=0$ in $\mbox{int}(\mathcal{B})\setminus\mathcal{A}$ and $z=0$ in $\mbox{int}(\mathcal{A})\setminus\mathcal{B}$. That is, $w\in H_0^1(\mbox{int}(\mathcal{A}))$ and $z\in H_0^1(\mbox{int}(\mathcal{B}))$.

On the one hand, one can get that $w_\e+z_\e$ satisfies, in the weak sense,
$$
-\Delta(w_\e+z_\e)\le \max\{\mu_1,\mu_2,\beta\}(w_\e+z_\e)^3,\,\ x\in\mathbb{R}^N.
$$
Recalling that $\{w_\e+z_\e\}$ is uniformly bounded in $H^1(\mathbb{R}^N)$ with respect to $\e$ and $\mu_1,\mu_2>0$ small, similarly as that in \cite[Proposition 2.2]{Barrios15}, one can get that  $\{w_\e+z_\e\}$ is uniformly bounded in $L^\infty(\Sigma)$ with respect to $\e$ and $\mu_1,\mu_2>0$ small. Moreover, $\|z\|_\infty+\|w\|_\infty\le C$ for some $C$(independent of $\mu_1,\mu_2$). On the other hand, denote by $\phi_1$ and $\lambda_1$ the first eigenfunction and eigenvalue of $(-\Delta, H_0^1(\mbox{int}(\mathcal{B})\setminus\mathcal{A})))$ respectively. Since
$$
-\Delta z_\e+\frac{b(x)}{\e^2}z_\e=\mu_2 z_\e^3+\beta w_\e^2 z_\e,\,\,x\in\mathbb{R}^N,
$$
we obtain that
$$
\lambda_1\int_{\mbox{int}(\mathcal{B})\setminus\mathcal{A}}z_\e \phi_1=\int_{\mbox{int}(\mathcal{B})\setminus\mathcal{A}))}(\mu_2 z_\e^3+\beta w_\e^2 z_\e)\phi_1.
$$
Taking the limit as $\varepsilon\rightarrow0$,
$$
\int_{\mbox{int}(\mathcal{B})\setminus\mathcal{A}}z \phi_1[\lambda_1-\mu_2 z^2]=0.
$$
So if $\mu_2>0$ is small such that $\lambda_1-\mu_2 z^2(x)>0$ a. e. in $\mbox{int}(\mathcal{B})\setminus\mathcal{A}$, we know $z=0$ in $\mbox{int}(\mathcal{B})\setminus\mathcal{A}$. That is, $z\in H_0^1(\mbox{int}(\Omega))$. Similarly, $w\in H_0^1(\mbox{int}(\Omega))$.

\vskip0.1in
{\bf Step 4.} For any $\eta,\varphi\in H_0^1(\mbox{int}(\Omega))$, we have
\begin{align*}
\int_{\mbox{int}(\Omega)}\nabla w_{\varepsilon}\nabla\eta+\nabla z_{\varepsilon}\nabla\varphi-\mu_{1}w_{\varepsilon}^{3}\eta-\beta z_{\varepsilon}^{2}w_{\varepsilon}\eta-\beta w_{\varepsilon}^{2}z_{\varepsilon}\varphi-\mu_{2}z_{\varepsilon}^{3}\varphi=J_\e'(w_\e,z_\e)(\eta,\varphi)=0
\end{align*}
and hence
\begin{align*}
\int_{\mbox{int}(\Omega)}\nabla w\nabla\eta+\nabla z\nabla\varphi-\mu_{1}w^{3}\eta-\beta z^{2}w\eta-\beta w^{2}z\varphi-\mu_{2}z^{3}\varphi=0.
\end{align*}
In view of the regularity assumptions on the set $\Omega$ and by classical regularity theory, $(w,z)$ is a fully nontrivial solution of problem \eqref{limitsystem} and
$J_\ast(w,z)\ge c_\ast$.

We also have
\begin{align*}
&\mathop{\limsup}\limits_{\varepsilon\rightarrow0}\int_{\mathbb{R}^{N}}|\nabla w_{\varepsilon}|^{2}+|\nabla z_{\varepsilon}|^{2}+w_\e^2+z_\e^2\\
&\leq\mathop{\lim}\limits_{\varepsilon\rightarrow0}\int_{\mathbb{R}^{N}}|\nabla w_{\varepsilon}|^{2}+\frac{a(x)}{\varepsilon^{2}}w_{\varepsilon}^{2}+|\nabla z_{\varepsilon}|^{2}+\frac{b(x)}{\varepsilon^{2}}z_{\varepsilon}^{2}+w_\e^2+z_\e^2\\
&=\mathop{\lim}\limits_{\varepsilon\rightarrow0}F(w_{\varepsilon},z_{\varepsilon})+\|w\|_2^2+\|z\|_2^2\\
&=F(w,z)+\|w\|_2^2+\|z\|_2^2=\int_{\mathbb{R}^{N}}|\nabla w|^{2}+|\nabla z|^{2}+w^2+z^2.
\end{align*}
This implies $(w_{\varepsilon},z_{\varepsilon})\rightarrow(w,z)$ strongly in $\mathbb{H}$ as $\varepsilon\rightarrow0$. Due to $d_\e=J_{\varepsilon}(w_\e,v_\e)\le c_\ast$, the following inequality holds
$$
J_\ast(w,z)\le\limsup_{\e\rightarrow0}J_\e(w_\e,z_\e)\le c_\ast,
$$
which yields that $(w,z)$ is a ground state of problem \eqref{limitsystem}. The proof is completed.
\ep


\end{document}